\input amstex
\documentstyle{amsppt}
\input epsf.tex

\def\stydate{May 10, 2002}

\chardef\tempcat\catcode`\@ \ifx\undefined\amstexloaded\input
amstex \else\catcode`\@\tempcat\fi \expandafter\ifx\csname
amsppt.sty\endcsname\relax\input amsppt.sty \fi
\let\tempcat\undefined

\immediate\write16{This is LABEL.DEF by A.Degtyarev <\stydate>}
\expandafter\ifx\csname label.def\endcsname\relax\else
  \message{[already loaded]}\endinput\fi
\expandafter\edef\csname label.def\endcsname{%
  \catcode`\noexpand\@\the\catcode`\@\edef\noexpand\styname{LABEL.DEF}%
  \def\expandafter\noexpand\csname label.def\endcsname{\stydate}%
    \toks0{}\toks2{}}
\catcode`\@11
\def\labelmesg@ {LABEL.DEF: }
{\edef\temp{\the\everyjob\W@{\labelmesg@<\stydate>}}
\global\everyjob\expandafter{\temp}}

\def\@car#1#2\@nil{#1}
\def\@cdr#1#2\@nil{#2}
\def\eat@bs{\expandafter\eat@\string}
\def\eat@ii#1#2{}
\def\eat@iii#1#2#3{}
\def\eat@iv#1#2#3#4{}
\def\@DO#1#2\@{\expandafter#1\csname\eat@bs#2\endcsname}
\def\@N#1\@{\csname\eat@bs#1\endcsname}
\def\@Nx{\@DO\noexpand}
\def\@Name#1\@{\if\@undefined#1\@\else\@N#1\@\fi}
\def\@Ndef{\@DO\def}
\def\@Ngdef{\global\@Ndef}
\def\@Nedef{\@DO\edef}
\def\@Nxdef{\global\@Nedef}
\def\@Nlet{\@DO\let}
\def\@undefined#1\@{\@DO\ifx#1\@\relax\true@\else\false@\fi}
\def\@@addto#1#2{{\toks@\expandafter{#1#2}\xdef#1{\the\toks@}}}
\def\@@addparm#1#2{{\toks@\expandafter{#1{##1}#2}%
    \edef#1{\gdef\noexpand#1####1{\the\toks@}}#1}}
\def\make@letter{\edef\t@mpcat{\catcode`\@\the\catcode`\@}\catcode`\@11 }
\def\donext@{\expandafter\egroup\next@}
\def\x@notempty#1{\expandafter\notempty\expandafter{#1}}
\def\lc@def#1#2{\edef#1{#2}%
    \lowercase\expandafter{\expandafter\edef\expandafter#1\expandafter{#1}}}
\newif\iffound@
\def\find@#1\in#2{\found@false
    \DNii@{\ifx\next\@nil\let\next\eat@\else\let\next\nextiv@\fi\next}%
    \edef\nextiii@{#1}\def\nextiv@##1,{%
    \edef\next{##1}\ifx\nextiii@\next\found@true\fi\FN@\nextii@}%
    \expandafter\nextiv@#2,\@nil}
{\let\head\relax\let\specialhead\relax\let\subhead\relax
\let\subsubhead\relax\let\proclaim\relax
\gdef\let@relax{\let\head\relax\let\specialhead\relax\let\subhead\relax
    \let\subsubhead\relax\let\proclaim\relax}}
\newskip\@savsk
\let\@ignorespaces\ignorespaces
\def\@ignorespacesp{\ifhmode
  \ifdim\lastskip>\z@\else\penalty\@M\hskip-1sp%
        \penalty\@M\hskip1sp \fi\fi\@ignorespaces}
\def\ignorespaces{\protect\@ignorespacesp}
\def\@bsphack{\relax\ifmmode\else\@savsk\lastskip
  \ifhmode\edef\@sf{\spacefactor\the\spacefactor}\fi\fi}
\def\@esphack{\relax
  \ifx\penalty@\penalty\else\penalty\@M\fi   
  \ifmmode\else\ifhmode\@sf{}\ifdim\@savsk>\z@\@ignorespacesp\fi\fi\fi}
\let\@frills@\identity@
\let\@txtopt@\identyty@
\newif\if@star
\newif\if@write\@writetrue
\def\@numopt@{\if@star\expandafter\eat@\fi}
\def\checkstar@#1{\DN@{\@writetrue
  \ifx\next*\DN@####1{\@startrue\checkstar@@{#1}}%
      \else\DN@{\@starfalse#1}\fi\next@}\FN@\next@}
\def\checkstar@@#1{\DN@{%
  \ifx\next*\DN@####1{\@writefalse#1}%
      \else\DN@{\@writetrue#1}\fi\next@}\FN@\next@}
\def\checkfrills@#1{\DN@{%
  \ifx\next\nofrills\DN@####1{#1}\def\@frills@####1{####1\nofrills}%
      \else\DN@{#1}\let\@frills@\identity@\fi\next@}\FN@\next@}
\def\checkbrack@#1{\DN@{%
    \ifx\next[\DN@[####1]{\def\@txtopt@########1{####1}#1}%
    \else\DN@{\let\@txtopt@\identity@#1}\fi\next@}\FN@\next@}
\def\check@therstyle#1#2{\bgroup\DN@{#1}\ifx\@txtopt@\identity@\else
        \DNii@##1\@therstyle{}\def\@therstyle{\DN@{#2}\nextii@}%
    \expandafter\expandafter\expandafter\nextii@\@txtopt@\@therstyle.\@therstyle
    \fi\donext@}

\newread\@inputcheck
\def\@input#1{\openin\@inputcheck #1 \ifeof\@inputcheck \W@
  {No file `#1'.}\else\closein\@inputcheck \relax\input #1 \fi}

\def\loadstyle#1{\edef\next{#1}%
    \DN@##1.##2\@nil{\if\notempty{##2}\else\def\next{##1.sty}\fi}%
    \expandafter\next@\next.\@nil\lc@def\next@\next
    \expandafter\ifx\csname\next@\endcsname\relax\input\next\fi}

\let\pagebody@\pagebody
\let\pagetop@\empty
\let\pagebot@\empty
\let\@Xend\empty
\def\pagebody{\pagetop@\pagebody@\pagebot@\@Xend}
\let\@Xclose\empty

\newwrite\@Xmain
\newwrite\@Xsub
\def\W@X{\write\@Xout}
\def\make@Xmain{\global\let\@Xout\@Xmain\global\let\end\endmain@
  \xdef\@Xname{\jobname}\xdef\@inputname{\jobname}}
\begingroup
\catcode`\(\the\catcode`\{\catcode`\{12
\catcode`\)\the\catcode`\}\catcode`\}12
\gdef\W@count#1((\lc@def\@tempa(#1)%
    \def\\##1(\W@X(\global##1\the##1))%
    \edef\@tempa(\W@X(%
        \string\expandafter\gdef\string\csname\space\@tempa\string\endcsname{)%
        \\\pageno\\\cnt@toc\\\cnt@idx\\\cnt@glo\\\footmarkcount@
        \@Xclose\W@X(}))\expandafter)\@tempa)
\endgroup
\def\readaux{\bgroup\checkbrack@\readaux@}
\let\begin\readaux
\def\readaux@{%
    \W@{>>> \labelmesg@ Run this file twice to get x-references right}%
    \global\everypar{}%
    {\def\\{\global\let}%
        \def\/##1##2{\gdef##1{\wrn@command##1##2}}%
        \disablepreambule@cs}%
    \make@Xmain{\make@letter\setboxz@h{\@input{\@txtopt@{\@Xname.aux}}%
            \lc@def\@tempa\jobname\@Name\open@\@tempa\@}}%
  \immediate\openout\@Xout\@Xname.aux%
    \immediate\W@X{\relax}\egroup}
\everypar{\global\everypar{}\readaux}
{\toks@\expandafter{\topmatter}
\global\edef\topmatter{\noexpand\readaux\the\toks@}}
\let\@@end@@\end

\def\@Xclose@{{\def\@Xend{\ifnum\insertpenalties=\z@
        \W@count{close@\@Xname}\closeout\@Xout\fi}%
    \vfill\supereject}}
\def\endmain@{\@Xclose@
    \W@{>>> \labelmesg@ Run this file twice to get x-references right}%
    \@@end@@}
\def\disablepreambule@cs{\\\disablepreambule@cs\relax}

\def\include#1{\bgroup
  \ifx\@Xout\@Xsub\DN@{\errmessage
        {\labelmesg@ Only one level of \string\include\space is supported}}%
    \else\edef\@tempb{#1}\clearpage
      \DN@##1 {\if\notempty{##1}\edef\@tempb{##1}\DN@####1\eat@ {}\fi\next@}%
    \DNii@##1.{\edef\@tempa{##1}\DN@####1\eat@.{}\next@}%
        \expandafter\next@\@tempb\eat@{} \eat@{} %
    \expandafter\nextii@\@tempb.\eat@.%
        \relaxnext@
      \if\x@notempty\@tempa
          \edef\nextii@{\write\@Xmain{%
            \noexpand\string\noexpand\@input{\@tempa.aux}}}\nextii@
        \ifx\undefined\@includelist\found@true\else
                    \find@\@tempa\in\@includelist\fi
            \iffound@\ifx\undefined\@noincllist\found@false\else
                    \find@\@tempb\in\@noincllist\fi\else\found@true\fi
            \iffound@\lc@def\@tempa\@tempa
                \if\@undefined\close@\@tempa\@\else\edef\next@{\@Nx\close@\@tempa\@}\fi
            \else\xdef\@Xname{\@tempa}\xdef\@inputname{\@tempb}%
                \W@count{open@\@Xname}\global\let\@Xout\@Xsub
            \openout\@Xout\@tempa.aux \W@X{\relax}%
            \DN@{\let\end\endinput\@input\@inputname
                    \@Xclose@\make@Xmain}\fi\fi\fi
  \donext@}
\def\includeonly#1{\edef\@includelist{#1}}
\def\noinclude#1{\edef\@noincllist{#1}}

\def\arabicnum#1{\number#1}

\def\Romannum#1{\expandafter\uppercase\expandafter{\romannumeral#1}}
\def\alphnum#1{\ifcase#1\or a\or b\or c\or d\else\@ialph{#1}\fi}
\def\@ialph#1{\ifcase#1\or \or \or \or \or e\or f\or g\or h\or i\or j\or
    k\or l\or m\or n\or o\or p\or q\or r\or s\or t\or u\or v\or w\or x\or y\or
    z\else\fi}
\def\Alphnum#1{\ifcase#1\or A\or B\or C\or D\else\@Ialph{#1}\fi}
\def\@Ialph#1{\ifcase#1\or \or \or \or \or E\or F\or G\or H\or I\or J\or
    K\or L\or M\or N\or O\or P\or Q\or R\or S\or T\or U\or V\or W\or X\or Y\or
    Z\else\fi}

\def\ST@P{step}
\def\ST@LE{style}
\def\N@M{no}
\def\F@NT{font@}
\outer\def\newcounter{\checkbrack@{\expandafter\newcounter@\@txtopt@{{}}}}
{\let\newcount\relax
\gdef\newcounter@#1#2#3{{%
    \toks@@\expandafter{\csname\eat@bs#2\N@M\endcsname}%
    \DN@{\alloc@0\count\countdef\insc@unt}%
    \ifx\@txtopt@\identity@\expandafter\next@\the\toks@@
        \else\if\notempty{#1}\global\@Nlet#2\N@M\@#1\fi\fi
    \@Nxdef\the\eat@bs#2\@{\if\@undefined\the\eat@bs#3\@\else
            \@Nx\the\eat@bs#3\@.\fi\noexpand\arabicnum\the\toks@@}%
  \@Nxdef#2\ST@P\@{}%
  \if\@undefined#3\ST@P\@\else
    \edef\next@{\noexpand\@@addto\@Nx#3\ST@P\@{%
             \global\@Nx#2\N@M\@\z@\@Nx#2\ST@P\@}}\next@\fi
    \expandafter\@@addto\expandafter\@Xclose\expandafter
        {\expandafter\\\the\toks@@}}}}
\outer\def\copycounter#1#2{%
    \@Nxdef#1\N@M\@{\@Nx#2\N@M\@}%
    \@Nxdef#1\ST@P\@{\@Nx#2\ST@P\@}%
    \@Nxdef\the\eat@bs#1\@{\@Nx\the\eat@bs#2\@}}
\outer\def\everystep{\checkstar@\everystep@}
\def\everystep@#1{\if@star\let\next@\gdef\else\let\next@\@@addto\fi
    \@DO\next@#1\ST@P\@}
\def\counterstyle#1{\@Ngdef\the\eat@bs#1\@}
\def\advancecounter#1#2{\@N#1\ST@P\@\global\advance\@N#1\N@M\@#2}
\def\setcounter#1#2{\@N#1\ST@P\@\global\@N#1\N@M\@#2}
\def\counter#1{\refstepcounter#1\printcounter#1}
\def\printcounter#1{\@N\the\eat@bs#1\@}
\def\refcounter#1{\xdef\@lastmark{\printcounter#1}}
\def\stepcounter#1{\advancecounter#1\@ne}
\def\refstepcounter#1{\stepcounter#1\refcounter#1}
\def\savecounter#1{\@Nedef#1@sav\@{\global\@N#1\N@M\@\the\@N#1\N@M\@}}
\def\restorecounter#1{\@Name#1@sav\@}

\def\warning#1#2{\W@{Warning: #1 on input line #2}}
\def\warning@#1{\warning{#1}{\the\inputlineno}}
\def\wrn@@Protect#1#2{\warning@{\string\Protect\string#1\space ignored}}
\def\wrn@@label#1#2{\warning{label `#1' multiply defined}{#2}}
\def\wrn@@ref#1#2{\warning@{label `#1' undefined}}
\def\wrn@@cite#1#2{\warning@{citation `#1' undefined}}
\def\wrn@@command#1#2{\warning@{Preamble command \string#1\space ignored}#2}
\def\wrn@@option#1#2{\warning@{Option \string#1\string#2\space is not supported}}
\def\wrn@@reference#1#2{\W@{Reference `#1' on input line \the\inputlineno}}
\def\wrn@@citation#1#2{\W@{Citation `#1' on input line \the\inputlineno}}
\let\wrn@reference\eat@ii
\let\wrn@citation\eat@ii
\def\nowarning#1{\if\@undefined\wrn@\eat@bs#1\@\wrn@option\nowarning#1\else
        \@Nlet\wrn@\eat@bs#1\@\eat@ii\fi}
\def\printwarning#1{\if\@undefined\wrn@@\eat@bs#1\@\wrn@option\printwarning#1\else
        \@Nlet\wrn@\eat@bs#1\expandafter\@\csname wrn@@\eat@bs#1\endcsname\fi}
\printwarning\Protect \printwarning\label \printwarning\ref
\printwarning\cite \printwarning\command \printwarning\option

{\catcode`\#=12\gdef\@lH{#}}
\def\@@HREF#1{}
\def\@HREF#1#2{\@@HREF{a #1}{\let\@@HREF\eat@#2}\@@HREF{/a}}
\def\@@Hf#1{file:#1} \let\@Hf\@@Hf
\def\@@Hl#1{\if\notempty{#1}\@lH#1\fi} \let\@Hl\@@Hl
\def\@@Hname#1{\@HREF{name="#1"}{}} \let\@Hname\@@Hname
\def\@@Href#1{\@HREF{href="#1"}} \let\@Href\@@Href
\ifx\undefined\pdfoutput
  \csname newcount\endcsname\pdfoutput
\else
  \def\pdflinkattr{attr{/C [0 0.9 0.9]}}
  \let\pdflinkbegin\empty
  \let\pdflinkend\empty
  \def\@pdfHf#1{file {#1}}
  \def\@pdfHl#1{name {#1}}
  \def\@pdfHname#1{\pdfdest name{#1}xyz\relax}
  \def\@pdfHref#1#2{\pdfstartlink \pdflinkattr goto #1\relax
    \pdflinkbegin#2\pdflinkend\pdfendlink}
  \def\@ifpdf#1#2{\ifnum\pdfoutput>\z@\expandafter#1\else\expandafter#2\fi}
  \def\@Hf{\@ifpdf\@pdfHf\@@Hf}
  \def\@Hl{\@ifpdf\@pdfHl\@@Hl}
  \def\@Hname{\@ifpdf\@pdfHname\@@Hname}
  \def\@Href{\@ifpdf\@pdfHref\@@Href}
\fi
\def\@Hr#1#2{\if\notempty{#1}\@Hf{#1}\fi\@Hl{#2}}
\def\@localHref#1{\@Href{\@Hr{}{#1}}}
\def\@countlast#1{\@N#1last\@}
\def\@@countref#1#2{\global\advance#2\@ne
  \@Nxdef#2last\@{\the#2}\@tocHname{#1\@countlast#2}}
\def\@countref#1{\@DO\@@countref#1@HR\@#1}

\def\Href@@#1{\@N\Href@-#1\@}
\def\Href@#1#2{\@N\Href@-#1\@{\@Hl{@#1-#2}}}
\def\Hname@#1{\@N\Hname@-#1\@}
\def\Hlast@#1{\@N\Hlast@-#1\@}
\def\cntref@#1{\global\@DO\advance\cnt@#1\@\@ne
  \@Nxdef\Hlast@-#1\@{\@DO\the\cnt@#1\@}\Hname@{#1}{@#1-\Hlast@{#1}}}
\def\HyperRefs#1{\global\@Nlet\Hlast@-#1\@\empty
  \global\@Nlet\Hname@-#1\@\@Hname
  \global\@Nlet\Href@-#1\@\@Href}
\def\NoHyperRefs#1{\global\@Nlet\Hlast@-#1\@\empty
  \global\@Nlet\Hname@-#1\@\eat@
  \global\@Nlet\Href@-#1\@\eat@}

\HyperRefs{label} {\catcode`\-11
\gdef\@labelref#1{\Hname@-label{r@-#1}}
\gdef\@xHref#1{\Href@-label{\@Hl{r@-#1}}} }
\HyperRefs{toc}
\def\@HR#1{\if\notempty{#1}\string\@HR{\Hlast@{toc}}{#1}\else{}\fi}



\def\bftext{\ifmmode\fam\bffam\else\bf\fi}
\let\@lastmark\empty
\let\@lastlabel\empty
\def\lastmark{\@lastmark}
\let\lastlabel\empty
\let\everylabel\relax
\let\everylabel@\eat@
\let\everyref\relax
\def\newlabel{\bgroup\everylabel\newlabel@}
\def\newlabel@#1#2#3{\if\@undefined\r@-#1\@\else\wrn@label{#1}{#3}\fi
  {\let\protect\noexpand\@Nxdef\r@-#1\@{#2}}\egroup}
\def\w@ref{\bgroup\everyref\w@@ref}
\def\w@@ref#1#2#3#4{%
  \if\@undefined\r@-#1\@{\bftext??}#2{#1}{}\else%
   \@xHref{#1}{\@DO{\expandafter\expandafter#3}\r@-#1\@\@nil}\fi
  #4{#1}{}\egroup}
\def\@@@xref#1{\w@ref{#1}\wrn@ref\@car\wrn@reference}
\def\@xref#1{\rom{\@@@xref{#1}}}
\let\xref\@xref
\def\pageref#1{\w@ref{#1}\wrn@ref\@cdr\wrn@reference}
\def\thepage{\ifnum\pageno<\z@\romannumeral-\pageno\else\number\pageno\fi}
\def\label@{\@bsphack\bgroup\everylabel\label@@}
\def\label@@#1#2{\everylabel@{{#1}{#2}}%
  \@labelref{#2}%
  \let\thepage\relax
  \def\protect{\noexpand\noexpand\noexpand}%
  \edef\@tempa{\edef\noexpand\@lastlabel{#1}%
    \W@X{\string\newlabel{#2}{{\@lastmark}{\thepage}}{\the\inputlineno}}}%
  \expandafter\egroup\@tempa\@esphack}
\def\label#1{\label@{#1}{#1}}
\def\fn@P@{\relaxnext@
    \DN@{\ifx[\next\DN@[####1]{}\else
        \ifx"\next\DN@"####1"{}\else\DN@{}\fi\fi\next@}%
    \FN@\next@}
\def\eat@fn#1{\ifx#1[\expandafter\eat@br\else
  \ifx#1"\expandafter\expandafter\expandafter\eat@qu\fi\fi}
\def\eat@br#1]#2{}
\def\eat@qu#1"#2{}
{\catcode`\~\active\lccode`\~`\@
\lowercase{\global\let\@@P@~\gdef~{\protect\@@P@}}}
\def\Protect@@#1{\def#1{\protect#1}}
\def\disable@special{\let\W@X@\eat@iii\let\label\eat@
    \def\footnotemark{\protect\fn@P@}%
  \let\footnotetext\eat@fn\let\footnote\eat@fn
    \let\refcounter\eat@\let\savecounter\eat@\let\restorecounter\eat@
    \let\advancecounter\eat@ii\let\setcounter\eat@ii
  \let\ifvmode\iffalse\Protect@@\@@@xref\Protect@@\pageref\Protect@@\nofrills
    \Protect@@\\\Protect@@~}
\let\notoctext\identity@
\def\W@X@#1#2#3{\@bsphack{\disable@special\let\notoctext\eat@
    \def\chapter{\protect\chapter@toc}\let\thepage\relax
    \def\protect{\noexpand\noexpand\noexpand}#1%
  \edef\next@{\if\@undefined#2\@\else\write#2{#3}\fi}\expandafter}\next@
    \@esphack}
\newcount\cnt@toc
\def\writeauxline#1#2#3{\W@X@{\cntref@{toc}\let\tocref\@HR}
  \@Xout{\string\@Xline{#1}{#2}{#3}{\thepage}}}
{\let\newwrite\relax
\gdef\@openin#1{\make@letter\@input{\jobname.#1}\t@mpcat}
\gdef\@openout#1{\global\expandafter\newwrite\csname
tf@-#1\endcsname
   \immediate\openout\@N\tf@-#1\@\jobname.#1\relax}}
\def\@@openout#1{\@openout{#1}%
  \@@addto\readaux@{\immediate\closeout\@N\tf@-#1\@}}
\def\auxlinedef#1{\@Ndef\do@-#1\@}
\def\@Xline#1{\if\@undefined\do@-#1\@\expandafter\eat@iii\else
    \@DO\expandafter\do@-#1\@\fi}
\def\beginW@{\bgroup\def\do##1{\catcode`##112 }\dospecials\do\@\do\"
    \catcode`\{\@ne\catcode`\}\tw@\immediate\write\@N}
\def\endW@toc#1#2#3{{\string\tocline{#1}{#2\string\page{#3}}}\egroup}
\def\do@tocline#1{%
    \if\@undefined\tf@-#1\@\expandafter\eat@iii\else
        \beginW@\tf@-#1\@\expandafter\endW@toc\fi
} \auxlinedef{toc}{\do@tocline{toc}}

\let\protect\empty
\def\Protect#1{\if\@undefined#1@P@\@\PROTECT#1\else\wrn@Protect#1\empty\fi}
\def\PROTECT#1{\@Nlet#1@P@\@#1\edef#1{\noexpand\protect\@Nx#1@P@\@}}
\def\pdef#1{\edef#1{\noexpand\protect\@Nx#1@P@\@}\@Ndef#1@P@\@}

\Protect\operatorname \Protect\operatornamewithlimits
\Protect\qopname@ \Protect\qopnamewl@ \Protect\text
\Protect\topsmash \Protect\botsmash \Protect\smash
\Protect\widetilde \Protect\widehat \Protect\thetag
\Protect\therosteritem
\Protect\Cal \Protect\Bbb \Protect\bold \Protect\slanted
\Protect\roman \Protect\italic \Protect\boldkey
\Protect\boldsymbol \Protect\frak \Protect\goth \Protect\dots
\Protect\cong \Protect\lbrace \let\{\lbrace \Protect\rbrace
\let\}\rbrace
\let\root@P@@\root \def\root@P@#1{\root@P@@#1\of}
\def\root#1\of{\protect\root@P@{#1}}

\def\frills{\ignorespaces\@txtopt@}
\def\frillsnotempty#1{\x@notempty{\@txtopt@{#1}}}
\def\numberline{\@numopt@}
\newif\if@theorem
\let\@therstyle\eat@
\def\@headtext@#1#2{{\disable@special\let\protect\noexpand
    \def\chapter{\protect\chapter@rh}%
    \edef\next@{\noexpand\@frills@\noexpand#1{#2}}\expandafter}\next@}
\let\AmSrighthead@\rightheadtext
\def\rightheadtext{\checkfrills@{\@headtext@\AmSrighthead@}}
\let\AmSlefthead@\leftheadtext
\def\leftheadtext{\checkfrills@{\@headtext@\AmSlefthead@}}
\def\@head@@#1#2#3#4#5{\@Name\pre\eat@bs#1\@\if@theorem\else
    \@frills@{\csname\expandafter\eat@iv\string#4\endcsname}\relax
        \ifx\protect\empty\@N#1\F@NT\@\fi\fi
    \@N#1\ST@LE\@{\counter#3}{#5}%
  \if@write\writeauxline{toc}{\eat@bs#1}{#2{\counter#3}\@HR{#5}}\fi
    \if@theorem\else\expandafter#4\fi
    \ifx#4\endhead\ifx\@txtopt@\identity@\else
        \headmark{\@N#1\ST@LE\@{\counter#3}{\frills\empty}}\fi\fi
    \@Name\post\eat@bs#1\@\ignorespaces}
\ifx\undefined\endhead\Invalid@\endhead\fi
\def\@head@#1{\checkstar@{\checkfrills@{\checkbrack@{\@head@@#1}}}}
\def\@thm@@#1#2#3{\@Name\pre\eat@bs#1\@
    \@frills@{\csname\expandafter\eat@iv\string#3\endcsname}
    {\@theoremtrue\check@therstyle{\@N#1\ST@LE\@}\frills
            {\counter#2}\@theoremfalse}%
    \@DO\envir@stack\end\eat@bs#1\@
    \@N#1\F@NT\@\@Name\post\eat@bs#1\@\ignorespaces}
\def\@thm@#1{\checkstar@{\checkfrills@{\checkbrack@{\@thm@@#1}}}}
\def\@capt@@#1#2#3#4#5\endcaption{\bgroup
    \edef\@tempb{\global\footmarkcount@\the\footmarkcount@
    \global\@N#2\N@M\@\the\@N#2\N@M\@}%
    \def\shortcaption##1{\global\def\sh@rtt@xt####1{##1}}\let\sh@rtt@xt\identity@
    \DN@{#4{\@tempb\@N#1\ST@LE\@{\counter#2}}}%
    \if\notempty{#5}\DNii@{\next@\@N#1\F@NT\@}\else\let\nextii@\next@\fi
    \nextii@#5\endcaption
  \if@write\writeauxline{#3}{\eat@bs#1}{{} \@HR{\@N#1\ST@LE\@{\counter#2}%
    \if\notempty{#5}.\enspace\fi\sh@rtt@xt{#5}}}\fi
  \global\let\sh@rtt@xt\undefined\egroup}
\def\@capt@#1{\checkstar@{\checkfrills@{\checkbrack@{\@capt@@#1}}}}
\let\captiontextfont@\empty

\ifx\undefined\subsubheadfont@\def\subsubheadfont@{\it}\fi
\ifx\undefined\proclaimfont\def\proclaimfont{\sl}\fi
\ifx\undefined\proclaimfont@\let\proclaimfont@\proclaimfont\fi
\def\proclaimfont{\proclaimfont@}
\ifx\undefined\definitionfont@\def\AmSdeffont@{\rm}
    \else\let\AmSdeffont@\definitionfont@\fi
\ifx\undefined\remarkfont@\def\remarkfont@{\rm}\fi

\def\newfont@def#1#2{\if\@undefined#1\F@NT\@
    \@Nxdef#1\F@NT\@{\@Nx.\expandafter\eat@iv\string#2\F@NT\@}\fi}
\def\newhead@#1#2#3#4{{%
    \gdef#1{\@therstyle\@therstyle\@head@{#1#2#3#4}}\newfont@def#1#4%
    \if\@undefined#1\ST@LE\@\@Ngdef#1\ST@LE\@{\headstyle}\fi
    \if\@undefined#2\@\gdef#2{\headtocstyle}\fi
  \@@addto\moretocdefs@{\\#1#1#4}}}
\outer\def\newhead#1{\checkbrack@{\expandafter\newhead@\expandafter
    #1\@txtopt@\headtocstyle}}
\outer\def\newtheorem#1#2#3#4{{%
    \gdef#2{\@thm@{#2#3#4}}\newfont@def#2#4%
    \@Nxdef\end\eat@bs#2\@{\noexpand\revert@envir
        \@Nx\end\eat@bs#2\@\noexpand#4}%
  \if\@undefined#2\ST@LE\@\@Ngdef#2\ST@LE\@{\proclaimstyle{#1}}\fi}}%
\outer\def\newcaption#1#2#3#4#5{{\let#2\relax
  \edef\@tempa{\gdef#2####1\@Nx\end\eat@bs#2\@}%
    \@tempa{\@capt@{#2#3{#4}#5}##1\endcaption}\newfont@def#2\endcaptiontext%
  \if\@undefined#2\ST@LE\@\@Ngdef#2\ST@LE\@{\captionstyle{#1}}\fi
  \@@addto\moretocdefs@{\\#2#2\endcaption}\newtoc{#4}}}
{
\outer\gdef\newtoc#1{{%
    \@DO\ifx\do@-#1\@\relax
    \global\auxlinedef{#1}{\do@tocline{#1}}{}%
    \@@addto\tocsections@{\make@toc{#1}{}}\fi}}}

\toks@\expandafter{\itembox@}
\toks@@{\bgroup\let\therosteritem\identity@\let\rm\empty
  \let\@Href\eat@\let\@Hname\eat@
  \edef\next@{\edef\noexpand\@lastmark{\therosteritem@}}\donext@}
\edef\itembox@{\the\toks@@\the\toks@}
\def\firstitem@false{\let\iffirstitem@\iffalse
    \global\let\lastlabel\@lastlabel}

\let\rosteritemrefform\therosteritem
\let\rosteritemrefseparator\empty
\def\rosteritemref#1{\hbox{\rosteritemrefform{\@@@xref{#1}}}}
\def\local#1{\label@\@lastlabel{\lastlabel-i#1}}

\def\xRef@P@{\gdef\lastlabel}
\def\xRef#1{\@xref{#1}\protect\xRef@P@{#1}}

\def\iref@P@{\gdef\lastref}
\def\itemref#1#2{\rosteritemref{#1-i#2}\protect\iref@P@{#1}}
\def\iref#1{\@xref{#1}\rosteritemrefseparator\itemref{#1}}

\def\eqref#1{\thetag{\@@@xref{#1}}}
\def\tagform@#1{\ifmmode\hbox{\rm\else\rom{\fi
        (\ignorespaces#1\unskip)\iftrue}\else}\fi}

\let\AmSfnote@\makefootnote@
\def\makefootnote@#1{\bgroup\let\footmarkform@\identity@
  \edef\next@{\edef\noexpand\@lastmark{#1}}\donext@\AmSfnote@{#1}}

\def\clearpage{\ifnum\insertpenalties>0\line{}\fi\vfill\supereject}

\def\proof{\checkfrills@{\checkbrack@{%
    \check@therstyle{\@frills@{\demo}{\frills{Proof}}{}}
        {\frills{}\envir@stack\endremark\envir@stack\enddemo}%
  \envir@stack\endproof\ignorespaces}}}
\def\everyendproof{\qed}
\def\endproof{\nofrillscheck{\frills@{\everyendproof}\revert@envir\endproof\enddemo}}

\let\AmSref\ref
\let\AmSrefstyle\refstyle
\let\plaincite\cite
\def\citei@#1,{\citeii@#1\eat@,}
\def\citeii@#1\eat@{\w@ref{#1}\wrn@cite\@car\wrn@citation}
\def\mcite@#1;{\plaincite{\citei@#1\eat@,\unskip}\mcite@i}
\def\mcite@i#1;{\DN@{#1}\ifx\next@\endmcite@
  \else, \plaincite{\citei@#1\eat@,\unskip}\expandafter\mcite@i\fi}
\def\endmcite@{\endmcite@}
\def\cite#1{\mcite@#1;\endmcite@;}
\PROTECT\cite
\def\refstyle#1{\AmSrefstyle{#1}\uppercase{%
    \ifx#1A\relax \def\@ref@##1{\AmSref\xdef\@lastmark{##1}\key##1}%
    \else\ifx#1C\relax \def\@ref@##1{\AmSref\no\counter\refno}%
        \else\def\@ref@{\AmSref}\fi\fi}}
\refstyle A
\newcounter\refno\null
\newif\ifRefs
\gdef\Refs{\checkstar@{\checkbrack@{\csname AmSRefs\endcsname
  \nofrills{\frills{References}%
  \if@write\writeauxline{toc}{vartocline}{\@HR{\frills{References}}}\fi}%
  \def\ref{\@ref@}\Refstrue\ignorespaces}}}
\let\ref\xref

\newif\iftoc
\pdef\tocbreak{\iftoc\hfil\break\fi}
\def\tocsections@{\make@toc{toc}{}}
\let\moretocdefs@\empty
\def\newtocline@#1#2#3{%
  \edef#1{\def\@Nx#2line\@####1{\@Nx.\expandafter\eat@iv
        \string#3\@####1\noexpand#3}}%
  \@Nedef\no\eat@bs#1\@{\let\@Nx#2line\@\noexpand\eat@}%
    \@N\no\eat@bs#1\@}
\def\MakeToc#1{\@@openout{#1}}
\def\newtocline#1#2#3{\Err@{\Invalid@@\string\newtocline}}
\def\make@toc#1#2{\penaltyandskip@{-200}\aboveheadskip
    \if\notempty{#2}
        \centerline{\headfont@\ignorespaces#2\unskip}\nobreak
    \vskip\belowheadskip \fi
    \@openin{#1}\relax
    \vskip\z@}
\def\contents{\readaux\checkfrills@{\checkbrack@{\@contents@}}}
\def\@contents@{\toc@{\frills{Contents}}\envir@stack\endcontents%
    \def\nopagenumbers{\let\page\eat@}\let\newtocline\newtocline@\toctrue
  \def\@HR{\Href@{toc}}%
  \def\tocline##1{\csname##1line\endcsname}
  \edef\caption##1\endcaption{\expandafter\noexpand
    \csname head\endcsname##1\noexpand\endhead}%
    \ifmonograph@\def\vartoclineline{\Chapterline}%
        \else\def\vartoclineline##1{\sectionline{{} ##1}}\fi
  \let\\\newtocline@\moretocdefs@
    \ifx\@frills@\identity@\def\\##1##2##3{##1}\moretocdefs@
        \else\let\tocsections@\relax\fi
    \def\\{\unskip\space\ignorespaces}\let\maketoc\make@toc}
\def\endcontents{\tocsections@\vskip-\lastskip\revert@envir\endcontents
    \endtoc}

\if\@undefined\selectf@nt\@\let\selectf@nt\identity@\fi
\def\Err@math#1{\Err@{Use \string#1\space only in text}}
\def\textonlyfont@#1#2{%
    \def#1{\RIfM@\Err@math#1\else\edef\f@ntsh@pe{\string#1}\selectf@nt#2\fi}%
    \PROTECT#1}
\tenpoint

\def\newshapeswitch#1#2{\gdef#1{\selectsh@pe#1#2}\PROTECT#1}
\def\shapeswitch#1#2#3{\@Ngdef#1\string#2\@{#3}}
\shapeswitch\rm\bf\bf \shapeswitch\rm\tt\tt
\shapeswitch\rm\smc\smc
\newshapeswitch\em\it
\shapeswitch\em\it\rm \shapeswitch\em\sl\rm
\def\selectsh@pe#1#2{\relax\if\@undefined#1\f@ntsh@pe\@#2\else
    \@N#1\f@ntsh@pe\@\fi}

\def\@itcorr@{\leavevmode
    \edef\prevskip@{\ifdim\lastskip=\z@ \else\hskip\the\lastskip\relax\fi}\unskip
    \edef\prevpenalty@{\ifnum\lastpenalty=\z@ \else
        \penalty\the\lastpenalty\relax\fi}\unpenalty
    \/\prevpenalty@\prevskip@}
\def\rom@P@#1{\@itcorr@{\selectsh@pe\rm\rm#1}}
\def\rom{\protect\rom@P@}
\def\Rom@P@#1{\@itcorr@{\rm#1}}
\def\Rom{\protect\Rom@P@}
{\catcode`\-11 \HyperRefs{idx} \HyperRefs{glo}
\newcount\cnt@idx \global\cnt@idx=10000
\newcount\cnt@glo \global\cnt@glo=10000
\gdef\writeindex#1{\W@X@{\cntref@{idx}}\tf@-idx
 {\string\indexentry{#1}{\Hlast@{idx}}{\thepage}}}
\gdef\writeglossary#1{\W@X@{\cntref@{glo}}\tf@-glo
 {\string\glossaryentry{#1}{\Hlast@{glo}}{\thepage}}}
}
\def\emph#1{\@itcorr@\bgroup\em\ignorespaces#1\unskip\egroup
  \DN@{\DN@{}\ifx\next.\else\ifx\next,\else\DN@{\/}\fi\fi\next@}\FN@\next@}
\def\makequoteactive{\catcode`\"\active}
{\makequoteactive\gdef"{\FN@\quote@}
\gdef\quote@{\ifx"\next\DN@"##1""{\quoteii{##1}}\else\DN@##1"{\quotei{##1}}\fi\next@}}
\let\quotei\eat@
\let\quoteii\eat@
\def\MakeIndex{\@openout{idx}}
\def\MakeGlossary{\@openout{glo}}

\def\endofpar#1{\ifmmode\ifinner\endofpar@{#1}\else\eqno{#1}\fi
    \else\leavevmode\endofpar@{#1}\fi}
\def\endofpar@#1{\unskip\penalty\z@\null\hfil\hbox{#1}\hfilneg\penalty\@M}

\newdimen\normalparindent\normalparindent\parindent
\def\firstparindent#1{\everypar\expandafter{\the\everypar
  \global\parindent\normalparindent\global\everypar{}}\parindent#1\relax}

\@@addto\disablepreambule@cs{%
    \\\readaux\relax
    \\\begin\relax
    \\\readaux@\relax
    \\\@openout\eat@
    \\\@@openout\eat@
    \/\Monograph\empty
    \/\MakeIndex\empty
    \/\MakeGlossary\empty
    \/\MakeToc\eat@
    \/\HyperRefs\eat@
    \/\NoHyperRefs\eat@
}

\csname label.def\endcsname


\def\punct#1#2{\if\notempty{#2}#1\fi}
\def\sppunct{\punct{.\enspace}}
\def\varpunct#1#2{\if\frillsnotempty{#2}#1\fi}

\def\headstyle#1#2{\numberline{#1\sppunct{#2}}\ignorespaces#2\unskip}
\def\headtocstyle#1#2{\numberline{#1\punct.{#2}}\space #2}

\def\specialtocstyle#1#2{#2}
\newcounter\section\null
\newcounter\subsection\section
\newcounter\subsubsection\subsection
\newhead\specialsection[\specialtocstyle]\null\endspecialhead
\newhead\section\section\endhead
\newhead\subsection\subsection\endsubhead
\newhead\subsubsection\subsubsection\endsubsubhead
\def\firstappendix{\global\sectionno0 %
  \counterstyle\section{\Alphnum\sectionno}%
    \global\let\firstappendix\empty}

\def\appendixtocstyle#1#2{\space\numberline{Appendix #1\sppunct{#2}}#2}
\newhead\appendix[\appendixtocstyle]\section\endhead

\let\endAmSdef\enddefinition
\def\proclaimstyle#1#2{\numberline{#2\varpunct{.\enspace}{#1}}\frills{#1}}
\copycounter\thm\subsubsection
\theorem\thm\endproclaim
\proposition\thm\endproclaim
\lemma\thm\endproclaim
\corollary\thm\endproclaim
\definition\thm\endAmSdef
\example\thm\endAmSdef

\def\captionstyle#1#2{\frills{#1}\numberline{\varpunct{ }{#1}#2}}
\newcounter\figure\null
\newcounter\table\null
\newcaption{Figure}\figure\figure{lof}\botcaption
\newcaption{Table}\table\table{lot}\topcaption

\copycounter\equation\subsubsection


\def\C{{\Bbb C}}
\def\R{{\Bbb R}}
\def\Z{{\Bbb Z}}
\def\Q{{\Bbb Q}}
\def\Rp#1{\R\roman P^{#1}}
\def\Cp#1{\C\roman P^{#1}}
\def\Im{\mathop{\roman{Im}}\nolimits}

\def\Ker{\mathop{\roman{Ker}}\nolimits}
\def\la{\langle}
\def\ra{\rangle}
\def\id{\mathop{\roman{id}}\nolimits}
\def\Aut{\mathop{\roman{Aut}}\nolimits}
\def\emptyset{\varnothing}
\def\oo{\varnothing}
\def\discr{\operatorname{discr}}
\def\Pic{\operatorname{Pic}}

\let\tm\proclaim
\let\endtm\endproclaim
\let\rk=\remark
\let\endrk=\endremark
\let\ge\geqslant
\let\le\leqslant
\let\+\sqcup
\let\dsum\+

\NoBlackBoxes \rightheadtext{Real four-dimensional cubic
hypersurfaces}

\topmatter
\title Deformation classes of real four-dimensional cubic hypersurfaces
\endtitle
\author S. Finashin, V. Kharlamov
\endauthor
\address Middle East Technical University,
Department of Mathematics\endgraf Ankara 06531 Turkey
\endaddress
\address
Universit\'{e} Louis Pasteur et IRMA (CNRS)\endgraf 7 rue Ren\'{e}
Descartes 67084 Strasbourg Cedex, France
\endaddress
\thanks
The second author is supported by ANR-05-BLAN-0053-01.
\endthanks
\subjclass 14P25, 14J10, 14N25, 14J35, 14J70
\endsubjclass
 \abstract We study real nonsingular cubic hypersurfaces $X\subset
P^5$ up to deformation equivalence combined with projective
equivalence and prove that they are classified by the conjugacy
classes of involutions induced by the complex conjugation in
$H_4(X)$. Moreover, we provide a graph $\Gamma_{K4}$ whose
vertices represent the equivalence classes of such cubics and
edges represent their adjacency. It turns out that the graph
$\Gamma_{K4}$ essentially coincides with the graph $\Gamma_{K3}$
characterizing a certain adjacency of real non-polarized
K3-surfaces.
 \endabstract
\endtopmatter

\rightline{\vbox{\hsize 70mm \noindent\eightit\baselineskip10pt
The most familiar logics in the modal family are constructed from
a weak logic called K (after Saul Kripke). Under the narrow
reading, modal logic concerns necessity and possibility. A variety
of different systems may be developed for such logics using K as a
foundation.
 \vskip3mm\noindent\eightrm
 Stanford Encyclopedia of Philosophy, Modal Logic. http://plato.stanford.edu/}}
 \vskip10mm

\document

\heading\S1. Introduction
\endheading

Hypersurfaces of degree $d$ in the projective space $P^n$ of
dimension $n$ naturally form a projective space $P^N$ of dimension
$N=\binom{n+d}d-1 $. This projective space has certain additional
structures that are due to its origin. In particular, it is
equipped with a special {\it discriminant hypersurface} $\Delta$
that consists of the polynomials $f\in P^N$ such that
$\{f=0\}\subset P^n$ is singular. In fact, $\Delta$ can be seen as
the locus of critical values of the projection $\Gamma\to P^N$ of
the universal (smooth) hypersurface
$$
\Gamma\subset P^N\times P^n,\quad \Gamma=\{\,(f,x)\,\vert\,
f(x)=0, f\in P^N, x\in P^n\}.
$$
As a consequence, over $P^N\setminus \Delta$ we get, in
Kodaira-Spencer terminology, a deformation family of nonsingular
varieties $\Gamma^0\to P^N\setminus\Delta$, where $\Gamma^0\subset
\Gamma$ is the preimage of $P^N\setminus\Delta$.

Everything said before works equally well over $\C$ and $\R$.
However, there is a principal difference: over $\C$ the space
$P^N\setminus \Delta$ is connected and all hypersurfaces $\{f=0\}$
with $f\in P^N\setminus \Delta$ are deformation equivalent, while
over $\R$ this is no longer the case (except $d=1$) and it is a
natural (and classical) task to understand the nature of the
connected components of $P^N(\R)\setminus \Delta(\R)$. We will
call two real nonsingular hypersurfaces in $P^n$ {\it deformation
equivalent} if they represent points in the same connected
component of $P^N(\R)\setminus \Delta(\R)$ (which is a natural
deformation equivalence for nonsingular hypersurfaces $X\subset
P^n$) and {\it coarse deformation equivalent} if one hypersurface
is deformation equivalent to a projective transformation of the
other (which is a natural deformation equivalence for embeddings
$X\to P^n$). In the latter case we speak on hypersurfaces in the
same {\it coarse deformation class}.

The group $PGL(n+1,\R)$ of real projective transformations of
$P^n$ is connected if $n$ is even, and has two connected
components otherwise. Therefore, a coarse deformation equivalence
coincides with the deformation equivalence if $n$ is even, and
otherwise a coarse deformation class contains at most two
deformation classes.

In this paper we treat the case of cubic hypersurfaces in $P^5$
and respond to the following questions: {\it how many coarse
deformation classes exist}\rom; and {\it how to distinguish
hypersurfaces of distinct coarse deformation classes} \rom?

Responses to similar questions for cubic hypersurfaces of lower
dimension go up to Newton (1704) in the case of curves, and to
Schl\"{a}fli and Klein (1858--1873) in the case of surfaces
(simple proofs going back to Newton \cite{New} and, respectively,
Klein \cite{Klein} can be found in \cite{DIK1}.). The case of
threefolds was treated recently by V.~Krasnov \cite{Kra}.

Two real nonsingular cubic curves in $P^2$ or surfaces in $P^3$
are deformation equivalent if and only if their real point sets
are homeomorphic. Even stronger, "homeomorphic" can be replaced by
"having the same homology groups". In fact, there are $2$ classes
of nonsingular cubic curves: $S^1$ and $S^1 \+S^1$, and $5$
classes of nonsingular cubic surfaces: $\Rp2 \+S^2$, $\Rp2$,
$\#_3\Rp2$, $\#_5\Rp2$, and $\#_7\Rp2$. Here, $\#$ stands for the
connected sum and $\+$, for the disjoint sum.

As is shown by V.~Krasnov, there are $9$ classes of real
nonsingular cubic threefolds in $P^4$, and they are distinguished
by the homology of the real point set plus one additional bit of
information: does the real point set realize a trivial or
non-trivial $\Z/2$-homology class in the complex point set. In
fact, his proof, similarly to Klein's, appeals to a relation
between the cubics in $P^{n}$ and complete intersections of
bi-degree $(2,3)$ in $P^{n-1}$ (in Krasnov's case $n=4$). Finally,
the classification of real nonsingular cubic threefolds happens to
be closely connected with the classification of real canonical
curves of genus $4$ in $P^3$, and to the fact that the latter
curves form $8$ deformation classes.

To solve the problem of coarse deformation classification of real
nonsingular cubic fourfolds in $P^5$, we need to appeal to the
complex point sets more systematically and somehow in a deeper
way.
 Namely, considering a real cubic fourfold $X$ we
associate with it the involution $c\: M\to M$ induced by the
complex conjugation on $M= H_4(X(\C))$. Let us call the
isomorphism class of this involution the {\it homological type} of
$X$. Our first goal is to prove that this invariant is sufficient
(in terminology of \cite{DIK2} it means that the cubic fourfolds
are homologically quasi-simple).

\tm{Theorem 1.1} If two real nonsingular cubic hypersurfaces in
$P^5$ have the same homological type, then they are coarse
deformation equivalent.
\endtm

In fact, we show that the homological type can be replaced by the
isomorphism class of the eigenlattice $M_-=\Ker (1+c)$. We also
enumerate all the possible eigenlattices,  which implies, in
particular, that the number of coarse deformation classes of real
nonsingular cubic hypersurfaces in $P^5$ is $75$.

In addition to a description of the coarse deformation classes of
real cubic fourfolds, our results provide the graph of their
adjacency. It turns out that this graph essentially coincides with
the adjacency graph for non-polarized real K3-surfaces, or equally
with the graph determined by real $6$-polarizations (in the latter
description, outgoing edges of the graph are real
$6$-polarizations of the real K3-surfaces representing the vertex;
see the precise definitions and formulations in section 3). Behind
these strong relations between deformation classes of cubic
fourfolds and K3-surfaces there are many reasons. One of them is a
certain close  relation between their lattices. To underline this
important similarity, we decided to call $K4$-{\it lattice} the
middle homology lattice of the cubic fourfolds.

The crucial role in our approach is played by Nikulin's coarse
deformation classification of $6$-polarized K3-surfaces in terms
of involutions on the K3-lattice and his results on the
arithmetics of such involutions, see \cite{N1,N2}.

The paper is organized as follows. In \S2, we discuss the central
projection correspondence for the nodal cubic hypersurfaces, which
yields the relation between the cubic fourfolds and $6$-polarized
K3-surfaces. The adjacency graphs of projective hypersurfaces and
the lattice graphs of involutions are introduced in \S3, where we
study the basic properties of these graphs, especially in the case
of cubic hypersurfaces in $P^5$ (K4-graph) and real involutions in
the K3-lattice (K3-graph). This section contains also a complete
description of the K3-graph. In \S4, we construct a morphism from
the K4-graph to the K3-graph, which allows us to reconstruct the
former from the latter, and therefore to prove the main theorem.
In conclusion, we give a few remarks, which contain in particular
some open questions. In the Appendix, we provide a list of the
real involutions in the K3-lattice and their eigenlattices; this
list is used in \S3 and \S4.

\subheading{Acknowledgements} The authors thank for the
hospitality IHES, where this research was initiated.

\heading \S2. The lattices of cubics
\endheading

\subheading{2.1. Cubics with a node or a cusp} The discriminant
hypersurface $\Delta\subset P^N$ can be naturally stratified. Its
principal stratum (the stratum of the highest dimension) is the
smooth part of $\Delta$, it is formed by those hypersurfaces which
have a non-degenerate double point (so called $A_1$ singular
point, or {\it node}) and no other singular points. In codimension
one there are two strata: one is formed by the hypersurfaces which
have two nodes and no other singular points; another is formed by
the hypersurfaces which have a {\it cusp} (quadratic suspension
over $x^3$, called also $A_2$) as their only singular point. We
denote by $\Delta_0$ the union of the principal and cuspidal
strata and by $\Delta_0(\R)$ its real part,
$\Delta_0(\R)=\Delta_0\cap\Delta(\R)$.

The following statement is well known.

\tm{Proposition 2.1} $\Delta_0(\R)$ is a topological submanifold
of $P^N(\R)$.\qed
\endtm

The connected components of $\Delta_0(\R)$ will be called {\it
pseudo-deformation classes} in $\Delta_0(\R)$ (the attribute
"pseudo" reflects a possibility of turning a node into a cusp).
The group $PGL(n+1,\R)$ naturally acts on $\Delta_0(\R)$ and we
call {\it projective classes} in $\Delta_0(\R)$ (in other words,
projective classes of cubic hypersurfaces with a node or cusp) the
orbits of this action. By {\it coarse pseudo-deformation classes}
in $\Delta_0(\R)$ we mean the orbits of the induced action on the
set of connected components of $\Delta_0(\R)$.

There is a natural {\it central projection correspondence},
between the projective classes of cubic hypersurfaces in $P^n$
with a marked double point and the projective classes of complete
intersections of bidegree $(2,3)$ in $P^{n-1}$.

Namely, start with  a cubic hypersurface $X$ in $P^n$ which has a
double point at $x\in X$, and consider an affine chart $A^n\subset
P^n$ centered at $x$. Then, the degree $3$ equation $f=0$ of $X$
written in such affine coordinates is defined up to a constant
factor and it splits, $f=f_2+f_3$, into two, quadratic and cubic,
homogeneous components $f_2,f_3$ defined up to a common factor. A
linear change of the affine coordinates leads to a simultaneous
projective transformation of $f_2, f_3$ (note that a homothety
with the center at $x$ multiplies $f_2, f_3$ by
$\lambda^2,\lambda^3$, where $\lambda\in\Bbb R$ in the case of
real hypersurfaces). A change of the affine chart centered at $x$
by a projective transformation preserving $x$ leads, in addition
to a projective transformation of $f_2, f_3$, to a transforming
$f_2, f_3$ into $f_2, f_3+l_1f_2$, where $1+l_1=0$ is a linear
equation defining the infinity hyperplane of the new affine chart.
Thus, the central projective correspondence yields a well defined
map which associates to the projective class of a cubic
hypersurface in $P^n$ with a marked double point a projective
class of complete intersections of bidegree $(2,3)$ in $P^{n-1}$.

 In general one should be careful about the difference between
the scheme-theoretic and the set-theoretic complete intersections
(which both appear in the above central projection
correspondence). For our purpose, {\it nonsingular complete
intersections} are sufficient. They can be defined as
set-theoretic intersections for which the Jacobian matrix of the
defining equations has the rank equal to the codimension. So
defined nonsingular complete intersections are scheme-theoretic
intersections: if $F_1,\dots,F_k$ are homogeneous polynomials in
$n$ variables such that the rank of their Jacobian matrix is equal
to $k$ at each point of $X\subset P^{n-1}$ defined by
$F_1=\dots=F_k=0$, then the ideal of $X$ is generated by
$F_1,\dots, F_k$.

\tm{Lemma 2.2}\label{nodeorcusp} The complete intersections which
are obtained by central projection from a cubic hypersurface $X$
in $P^n$ with a marked double point are nonsingular complete
intersections if and only if $X\in \Delta_0$.
\endtm

\demo{Proof} (cf., \cite{Kra}) Consider an affine chart
$A^n\subset P^n$ centered at the double point and denote by
$f=f_2+f_3$ the degree 3 polynomial defining the cubic. A
straightforward check shows that the singularity at $0\in A^n$ of
$f_2+f_3=0$ is not of the type $A_1$ or $A_2$ if $f_2=f_3=0$ has a
singular point which is also a singular point of $f_2=0$ (that is
a point with $df_2=f_2=f_3=0$), and that $f_2+f_3=0$ has more than
one singular point if $f_2=f_3=0$ has a singular point which is a
nonsingular point of $f_2=0$ (that is a point with  $df_2\ne 0$).
Finally, if $f_2+f_3=0$ has a singular point $x$ different from
$0$, then $df_2+df_3=0$ at $x$, hence, by the Euler relation,
$2f_2+3f_3=0$ at $x$, which implies $f_2=f_3=0$ at $x$, so that
the intersection $f_2=f_3=0$ is singular.
 \qed\enddemo

The group $PGL(n,\R)$ naturally acts on the set of real
nonsingular complete intersections of bi-degree $(2,3)$ in
$P^{n-1}$ and we call {\it projective classes of real nonsingular
complete intersections of bi-degree $(2,3)$} the orbits of this
action.

\tm{Lemma 2.3}\label{projclasses} The central projection
correspondence defines a bijection between the set of projective
classes in $\Delta_0(\R)$ and the set of projective classes of
real nonsingular complete intersections of bi-degree $(2,3)$.
\endtm

\demo{Proof} Performing a projective transformation of $f_2,f_3$
leads to a projective transformation of $f=f_2+f_3$. Since
$f_2,f_3$ generate the ideal of the intersection, to obtain
bijectivity it remains to check the effect of independent
rescaling of $f_2,f_3$ and also of replacing $f_2,f_3$ by $f_2,
f_3+l_1f_2$, where $l_1$ is a degree $1$ homogeneous polynomial.
As we have seen above analyzing the central correspondence, the
both modifications lead to a combination of rescaling and a
projective transformation of $f=f_2+f_3$.\qed
\enddemo

The central projection correspondence extends to real algebraic
(or real analytic) families of cubics with a marked double point,
where by a real algebraic (respectively, real analytic) family of
cubics with a marked double point we mean a family which can be
given by a real algebraic (respectively, real analytic) family of
cubic equations equipped with a real algebraic (respectively, real
analytic) family of their double points. Namely, a local real
family of cubics can be calibrated by a real family of projective
transformations  to have the marked double points at the same
point, and then the family of cubics can be defined by a real
family of equations $f_2+f_3$ so that we get in $P^{n-1}$ a well
defined, up to projective equivalence, real family of complete
intersections, $f_2=f_3=0$. In particular, in accordance with
Lemma 2.2, there is a well defined map which associates to the
coarse pseudo-deformation class of a real cubic hypersurface in
$P^n$ with a node or cusp a coarse deformation class of real
nonsingular complete intersections of bidegree $(2,3)$ in
$P^{n-1}$ (here, as usual, we call two real nonsingular complete
intersections coarse deformation equivalent if one of them can be
connected with a projective transformation of another by a real
family of nonsingular complete intersections).

\tm{Proposition 2.4} The coarse pseudo-deformation classes in
$\Delta_0(\R)$ are in one-to-one correspondence with the coarse
deformation classes of nonsingular complete intersections of
bi-degree $(2,3)$ in $P^{n-1}(\R)$.
\endtm

\demo{Proof} Any local real deformation family of nonsingular
complete intersections of bi-degree $(2,3)$ in $P^{n-1}$ can be
lifted to a family of (degree 2 and 3) homogeneous polynomials,
$f_2$ and $f_3$, generating the ideal, and thus to a real family
of cubics $f_2+f_3$, which according to Lemma 2.2 have a node or a
cusp as their unique singular point. Thus, the result follows from
Lemma 2.3. \qed
\enddemo

\subheading{2.2. Homology of singular cubics} Suppose that a cubic
$X_0$ has a node or a cusp at $x\in X_0$ and that it is the only
singular point. The central projection of $X_0$ from $x$ to
$P^{n-1}$ is a birational isomorphism. The singularity of $X_0$ at
$x$ is resolved by a simple blowup
 of $P^{n}$ at $x$, and the
central projection $X_0\to P^{n-1}$ lifts to a regular map $\pi:
\widehat X_0\to P^{n-1}$ which is the blowup of $P^{n-1}$ at the
(codimension two, nonsingular) complete intersection $Y$
representing $X_0$.

The next proposition applied to $B=P^{n-1}$ and $X=\hat X_0$ shows
how are related the middle dimensional lattices of $\hat X_0$ and
$Y$.

In this proposition, and throughout the paper, by {\it a real
structure} on a complex (algebraic or analytic) variety $B$ we
mean an antiholomorphic involution $B\to B$.

\tm{Proposition 2.5} Let $p\:X\to B$ denote a blowup of a compact
complex variety $B$ along a non-singular subvariety $Y\subset B$.
Then there exists a short exact sequence
$$
\minCDarrowwidth{3mm} \CD
0@>>>H^*(B)@>p^*\oplus{\operatorname{in}^*}>>H^*(X)\oplus
H^*(Y)@>\operatorname{in}^*+p^*\vert_Y>>H^*(p^{-1}(Y))@>>>0\endCD
$$
 If $Y$ is of codimension $2$ in
$B$, this exact sequence implies that
$$
\CD H^*(B)\oplus
H^{*}(Y)@>p^*+\operatorname{in}^!p^*\vert_Y>>H^*(X)\endCD
$$
is an isomorphism, and, in particular, for any $k$ the group
$H_k(X)$ splits as $H_{k-2}(Y)\oplus H_k(B)$. In the middle
dimension, for $\dim B=2m$, it gives a lattice isomorphism between
$H_{2m}(X)$ and $-(H_{2m-2}(Y))\oplus H_{2m}(B)$.

If $Y$ is a real non-singular codimension $2$ subvariety of a
complex non-singular variety $B$ with a real structure $c_B:B\to
B$, then  the induced involutions $c_X$ and $c_Y$ on the lattices
are related as $c_X(x\oplus y)=-c_Y(x)\oplus c_B(y)$, for any
$x\in H_{2m-2}(Y)$, $y\in H_{2m}(B)$.
\endtm

\demo{Proof}(cf. \cite{A}) The short exact sequence follows from
the commutativity of the diagram of long exact sequences of the
pairs $(B,Y)$ and $(X,p^{-1}(Y))$ due to the injectivity of the
pull back $H^*(B)\to H^*(X)$. To deduce the splitting $H^*(B)
\oplus H^{*}(Y)= H^*(X)$ it is sufficient to use the Leray-Hirsch
description of $H^*(p^{-1}(Y))$. The manifold $p^{-1}(Y)$ is the
projectivization $\Bbb P(\Cal N)$ of the normal bundle $\Cal N$ of
$Y$ in $B$, and the normal bundle of $p^{-1}(Y)$ in $X$ is nothing
but $O_{\Bbb P(\Cal N)}(-1)$ (that is the dual of the tautological
line bundle of $\Bbb P(\Cal N)$). All these constructions respect
the real structure.\qed
\enddemo

As it follows for example from Propositions 2.4 and 2.5, all
cubics $X\in \Delta_0(\R)$ from the same pseudo-deformation class
yield the same real homology type of $\hat X$ and in particular,
the same middle dimensional lattice and the same action of the
real structure involution on the lattice.

\subheading{2.3. Lattice twists} Given an integral quadratic form
$q\:L\to\Z$ on a finite rank free abelian group $L$, and $v\in L$
with $q(v)=\pm2$, there is a unique integral quadratic form $q'
\:L\to\Z$ such that $q'(v)=-q(v)$ and $q'(w)=q(w)$ for all $w\in
v^\perp$. Such a new quadratic form is given by $q'(x)=b(x,s_vx)$,
where $b$ is the bilinear form associated with $q$ and $s_v$ is
the reflection $L\to L$ defined by $x\mapsto x-b(x,v)v$ if
$q(v)=2$ and by $x\mapsto x+b(x,v)v$ if $q(v)=-2$. We call $q'$
{\it the $v$-twist of $q$} and denote by $t_v(L)$ the lattice
$(L,q')$.

Clearly, each automorphism $g:L\to L$ of $q$ with $gv=\pm v$ is an
automorphism of $q'$ (and vice versa).

\tm{Lemma 2.6} The quadratic forms $q'$ and $q$ have the same
discriminant group. In particular, $q'$ is non-degenerate
(unimodular) if, and only if, $q$ is non-degenerate (unimodular).
\endtm

\demo{Proof} According to $q'(x)=b(x,s_vx)$, the correlation
homomorphisms $L\to L^*$ of $q$ and $q'$ are related by the
reflection $s_v$ and, hence, have the same image. It remains to
notice that the discriminant group is the cokernel of the
correlation homomorphism.\qed
\enddemo

In what follows we abbreviate $q(x)=x^2$ and $b(x,y)= \la
x,y\ra=xy$ when it does not lead to a confusion.We denote the
discriminant group of $L$ by $\discr L$ or $\frak L$. Recall that
$\frak L$ is a finite group, if $q$ is non-degenerate, and then
$b$ induces a finite symmetric bilinear form $\frak b\:\frak
L\times\frak L\to\Q/\Z$. If, in addition, $L$ is even (that is
$x^2=0\mod 2$ for any $x\in L$), the discriminant group carries a
canonical finite quadratic form $\frak q:\frak L\to\Q/2\Z$ defined
via $\frak q(x+L)=x^2\mod 2\Z$. To extend this construction to odd
lattices it is sufficient to select a characteristic element $w\in
L\otimes\Q$ (that is an element $w$ such that $b(w, l)=l^2\mod
2\Z$ for any $l\in L$) and to put $\frak q(x+L)=x^2+b(x,w)\mod
2\Z$.

A finite rank free abelian group $L$ equipped with a symmetric
bilinear form $q$ is called {\it lattice}. If the form is
non-degenerate, we say that the lattice is non-degenerate.

\subheading{2.4. Homology of non-singular cubics} Our next goal is
to relate the middle dimensional lattice of $\widehat X_0$ with
that of a non-singular perturbation $X$ of $X_0$. This relation
holds for any hypersurface $X_0\subset P^{2m+1}$ having a single
node or a cusp. Since all the hypersurfaces from the same
connected component of $\Delta_0(\R)$ have the same middle
dimensional lattice of $\hat X_0$ (and the same cohomology ring)
with the same action of the real structure, we will restrict
ourselves to the case of hypersurfaces with a node.

Recall that given a generic one-parameter perturbation, $t\mapsto
X_t\in P^N$, $t\in\C$, $\vert t\vert<\epsilon$, of a hypersurface
$X_0$ with a node, there is a well defined up to sign vanishing
class $v_t\in H_{2m}(X_t)$ with $(-1)^m v_t^2=2$ (this class is
realized by so called {\it vanishing spheres}, which are
lagrangian spheres in $X_t$, and, hence, their normal bundle is
isomorphic to the cotangent bundle). If $X_0$ is real and belongs
to $\Delta_0(\R)$, for any real $t\ne 0$, then $c_t(v_t)=\pm v_t$,
where the sign $\pm$ depends only on the side of $\Delta_0(\R)$,
to where $X_0$ is shifted by the perturbation $X_t$. Moreover,
there are canonical, well defined up to isotopy, diffeomorphisms
between non-singular fibers $X_t$ with $\Im t\ge 0$, $t\ne 0$,
such that for any $t>0$ it holds $c_t\circ c_{-t}=\mu$, where
$\mu$ is the one-turn monodromy (also well defined up to isotopy).
These diffeomorphisms provide an identification of the
Poincar\'{e} duality lattices $M=H_{2m}(X_t), t\ne 0, t\in\R$,
together with the vanishing cycles $v=v_t$,  and in the case of a
single node, the monodromy $\mu$ is the Picard-Lefschetz
transformation $s_v:x\mapsto x-(-1)^m\la x,v\ra v$. In particular,
in this case $c_+=(c_{t>0})_* :M\to M$ commutes with
$c_-=(c_{t<0})_*:M\to M$ and they coincide on the orthogonal
complement $v^\perp$ of $v$, as it follows from $c_+\circ
c_-=s_v$. The vanishing cycle $v$ belongs to one of the
eigenlattices $M_\pm(c_+)=\Ker(1\mp c_+)$ and jumps to the
complementary sublattice $M_\mp(c_{-})=\Ker(1\pm c_-)$ of the
other involution.

The comparison of the Euler characteristic of $X_t(\R)$ with that
of $X_{-t}(\R)$ allows us to distinguish the two sides of
$D_0(\R)$, or in the other words, to coorient $D_0(\R)$. Such a
coorientation can be also characterized by the action of $c_t$ on
the vanishing class, as we formulate in the following Proposition.

\tm{Proposition 2.7} Given a real $t$, $0<\vert t\vert<\epsilon $,
the equality $c_t(v_t)=v_t$ holds if and only if
$\chi(X_t(\R))>\chi(X_{-t}(\R))$ {\rm(}which is equivalent to
$\chi(X_{-t}(\R))=\chi(\hat X_0(\R))$\rm ).
\endtm

\demo{Proof} All the statements follow from the direct inspection
of the {\it local Euler characteristic} of the real loci. By the
latter we mean the Euler characteristic of the real part of the
Milnor fiber, $X_t^{loc}(\R)\subset X_t(\R)$, in the case of $t\ne
0, t\in\R$, and the Euler characteristic of the real part of a
regular neighborhood of the exceptional divisor, $\hat
X_0^{loc}(\R)\subset \hat X_0(\R)$, in the case of $t=0$. The
Morse lemma implies that
$\chi(X_t(\R))-\chi(X_{-t}(\R))=\chi(X_t^{loc}(\R))-\chi(X_{-t}^{loc}(\R))=2$
if  $c_t(v_t)=v_t$. \qed
\enddemo

\tm{Corollary 2.8} Increase of $\chi(X_t(\R))$ defines a
co-orientation of $D_0(\R)$. For the increasing Euler
characteristic perturbations $X_t$, $t\ne 0$, it holds $\chi(\R
X_t)=\chi(\R X_{-t})+2$. \qed \endtm

It will be convenient for us to use such an ``increasing''
coorientation in case of odd $m$ and the opposite, ``decreasing''
coorientation in case of even $m$. The side of $D_0(\R)$ and the
corresponding perturbation $X_t$ for which $c_t(v_t)=(-1)^mv_t$
will be called {\it the ascendant side} and respectively {\it the
ascendant perturbation}.

Recall next that the exceptional divisor
$Q=p^{-1}(x_0)\subset\widehat X_0$ is a $(2m-1)$-dimensional
quadric in $P^{2m}$ and the normal bundle of $Q$ in $\widehat X_0$
is isomorphic to the line bundle induced from $O_{P^{2m}}(-1)$.
Thus, the generator $w_Q\in H_{2m}(Q)$ realized by intersection of
$Q$ with a codimension $m-1$ plane in $P^{2m}$ provides an element
$w=\text{in}_*w_Q\in H_{2m}(\widehat X_0)$ with $w^2=-2$.

\tm{Proposition 2.9} If $m$ is even, then the lattices
$M=H_{2m}(X_t), t\ne 0,$ and $\widehat M=H_{2m}(\widehat X_0)$ are
related by a $v$-twist about a vanishing cycle $v\in M$.
\endtm

\demo{Proof} It is sufficient to compare the exact sequences
$$
\align H^{2m-1}(V)=&0\to H^{2m}(X_t,V)\to H^{2m}(X_t)\to
H^{2m}(V)\cong\Z, \\
H^{2m-1}(Q)=&0\to H^{2m}(\widehat X_0,Q)\to H^{2m}(\widehat
X_0)\to H^{2m}(Q)\cong\Z,
\endalign
$$
where $V$ stands for a vanishing sphere, and to observe that the
both homomorphisms $H^{2m}(V)\to H^{2m+1}(X_t,V)$ and
$H^{2m}(Q)\to H^{2m+1}(\widehat X_0,Q)$ vanish, since $v=[V]\in
H_{2m}(X_t)$ and $w\in H_{2m}(\widehat X_0)$ are primitive (which
follows from $v^2= 2$ and $w^2=-2$).\qed
\enddemo

In what follows we fix a vanishing class $v$ (which is
well-defined only up to sign) and select that of the two group
isomorphisms $M=\widehat M$ given by Proposition 2.9 for which
$v=w$.

\tm{Corollary 2.10} Let $X_0$ be a cubic with a single node,
$Y\subset P^{2m}$ be the complete intersection corresponding to
$X_0$, and $X=X_t$ with $t\ne 0$. If $m$ is even, then the
lattices of $X, \hat X_0, Y$, and $P^{2m}$ are related as follows:
$$\align
H_{2m}(\widehat X_0)&=(-H_{2m-2}(Y))\oplus H_{2m}(P^{2m}(\C)),\\ 
H_{2m}(X)&=t_w(H_{2m}(\widehat X_0)),\quad w=h+2e,
\endalign
$$
where $h\in H_{2m-2}(Y)$ is the hyperplane-section class  and $e$
the canonical generator $[ P^{m}(\C)]\in H_{2m} (P^{2m}(\C))=\Z$.
All the relations respect the real structure, and, in particular,
if $X$ is obtained by an ascendant real perturbation of $\hat X_0$
then the involutions $c_X : H_{2m}(X)\to H_{2m}(X)$ and $c_{\hat
X_0}: H_{2m}(\widehat X_0)\to H_{2m}(\widehat X_0)$ induced by the
real structures coincide as we identify $H_{2m}(X)$ and
$H_{2m}(\widehat X_0)$ as groups.
\endtm

\demo{Proof} All the statements, except the formula $w=h+2e$,
follow directly from Propositions 2.5 and 2.9. This formula
follows from $H_{2m}(\widehat X)=(-H_{2m-2}(Y))\oplus
H_{2m}(P^{2m}(\C))$, since $w=\text{in}_*w_Q\in H_{2m}(\widehat
X)$, where $w_Q\in H_{2m}(Q)$ is realized by intersection of $Q$
with a codimension $m-1$ plane in $P^{2m}$.\qed
\enddemo

\subheading{2.5. The lattices K3 and K4} Let $X_0$ be a cubic
hypersurface in $P^5$ with a unique singular point which is a
node, and $Y\subset P^4$ be associated with $X_0$ (so that $m=2$
in notation of the previous sections). Denote by $\hat M$ and $M$
the lattices $\hat M=H_4(\hat X)$ and $M=H_4(X)$, where $X$ is a
perturbation of $X_0$.

According to Lemma 2.2, $Y$ is non-singular. Therefore, $Y$ is a
$K3$-surface which is $6$-polarized by a very ample hyperplane
section divisor $P, P^2=6$ (it may be worth to mention here that
conversely any K3-surface of degree 6 in $P^4$ is a complete
intersection of a quadric and a cubic, see, e.g., \cite{StD}). In
particular, $L=H_2(Y)$ is the K3-lattice $3U\oplus 2E_8$ ($U$
states for the unique even unimodular lattice of rang $2$ (and
signature $(1,1)$), and $E_8$ for the unique even unimodular
negative definite lattice of rang $8$) with a marked element
$h=[P]$, $h^2=6$.

By Corollary 2.10,
$$\align
\hat M&= (-L)\oplus\Z, \\
M&=t_w(\hat M), \quad w=h+2e,
\endalign$$
where $e$ is the generator of $\Z$ (corresponding to $[\Cp2]\in
H_4(\Cp4)$) and $e^2=1$. Note that $h^2=-6$ in $\hat M$.

If $X_0$ and $X$ are real, then $\hat X_0$ and $Y$ are real as
well. By Corollary 2.10, the action of the real structures on
$\hat M$ and $L$ are related as $c_{\hat M}(x\oplus
ne)=-c_L(x)\oplus ne$, and if the vanishing class $v\in M$ of a
real deformation $X$ is $c_M$-invariant (the case of ascendant
perturbations), then the involutions $c_M(x)$ and $c_{\hat M}(x)$
coincide as we identify $M$ and $\hat M$ as groups.

Let us denote by $\hat M_\pm$, $M_\pm$ and $L_\pm$ the $(\pm
1)$-eigenlattices of the involutions induced on $\hat M$, $M$ and
$L$ by the real structures. The following relations are also
immediate consequences of Corollary 2.10.

 \tm{Lemma 2.11}
If $X_0$ is real, then
$$\align
\hat M_+=&(-L_-)\oplus\Z.\\
\hat M_-=&-L_+,\\
\endalign
$$
and if, an addition, $M$ corresponds to an ascendant real
perturbation $X$ of $X_0$, then
$$
\align
M_+=&t_w(-L_-\oplus\Z),\\
M_-=&-L_+\,.  \qed
\endalign
$$
\endtm

The lattices $\hat M$ and $M$ are $3$-polarized by the fundamental
class $H$ of the intersection of $\hat X_0$ and $X$ respectively
with a generic projective subspace of codimension $2$ (in other
words, by the square of the hyperplane section).

\tm{Lemma 2.12} The vanishing class $v\in M$ and polarization
class $H\in M$ can be expressed as $v= h+2e$, $H=h+3e$ under the
group identification $M=\hat M=L\oplus \Z$.
\endtm

\demo{Proof} The first relation follows from $v=w$ and $w=h+2e$.

To prove the other one it is sufficient to check that $\la
H,\xi\ra=\la h+3e,\xi\ra$ for $\xi=e$ and also for any $\xi\in
-L$. Note that under identification $\hat M=M$ the polarization
class $H\in M$ becomes $D\circ D$, where $D\in H_6(\hat X_0)$ is
the fundamental class of the proper transform $C\subset \hat X_0$
of the cubic threefold given in $P^4$ by the equation $f_3=0$.
Since $e$ is realized by a generic $2$-plane in $P^4$, we have
$\la H, e\ra=3=\la h+3e, e\ra$. If $\xi\in -L$, then $\la
e,\xi\ra=0$ and under the group identification $-L=H_2(Y)$ we have
$\la h+3e,\xi\ra=\la -H_Y,\xi\ra$, where $H_Y\in H_2(Y)$ is the
fundamental class of a hyperplane section of $Y$. Therefore, it
remains to notice that the normal bundle of $C\cap p^{-1}(Y)$ in $
p^{-1}(Y)=\text{Proj} (O(3)\oplus O(2))$ is $O(-1)$, so that $\la
H,\xi\ra= \la D\circ D, p^*\xi\ra =\la -H_Y,\xi\ra$.\qed
\enddemo

Finally, let us observe that the both lattices $M$ and $\hat M$
are odd, with the natural (Chern) representatives of the
characteristic classes belonging to $M_+$ and, respectively, to
$\hat M_+$, so that the lattices $M_-$ and $\hat M_-$ are even.

\subheading{2.6. Bi-nodal $4$-cubics and flips of 6-polarized
$K3$} Assume that a cubic $X$ in $P^5$ has two non-degenerate
nodes, $s_1$ and $s_2$, and that it has no other singular points.
Pick projective coordinates $(x_0,x_1,\dots,x_5)$ such that
$s_1=(1,0,0,\dots,0)$ and $s_2=(0,1,0,\dots,0)$. Then, in the
affine coordinates centered at $s_1$ (here we let $x_0=1$) $X$ is
given by equation
$$
x_1L+A+(x_1B+C),
$$
and in the affine coordinates centered at $s_2$ (here we let
$x_1=1$) it is given by equation
$$
x_0L+B+(x_0A+C),
$$
where $L,A,B,$ and $ C$ depend only on $x_2,\dots, x_5$ and have
degree $1,2,2,$ and $3$ respectively.

As usual, we associate with a node $s_i$ of $X$ a $6$-polarized
K3-surface $Y_i$ which is given in $P^4$ by system of equations
$x_1L+A=x_1B+C=0$, if $i=1$, and $x_0L+B=x_0A+C=0$, if $i=2$. Each
of $Y_i$ has a single singular point, it is a nondegenerate double
point and is located at the center of the affine chart. The
surfaces $Y_i$ are birationally equivalent and a birational
isomorphism between them is given by $x_1=x_0\frac{A}B$.

Consider the induced biregular isomorphism $f:\hat Y_1\to\hat Y_2$
between the minimal models of $Y_1,Y_2$ and note that as any
birational isomorphism between minimal $K3$-surfaces $f$ is
biregular. In our case, each $\hat Y_i\to Y_i$ contracts a unique
$(-2)$-curve. Denote their homology classes by $v_i$, and the
hyperplane section homology classes by $h_i$.

\tm{Lemma 2.13}
$$\align
h_2&=f_*(2h_1-3v_1),\\
v_2&=f_*(h_1-2v_1).
\endalign
$$
\endtm

\demo{Proof} The classes $h_i-v_i$ are given by the proper
transform of hyperplane sections passing through $s_i$. Since $f$
preserves this condition, it implies $f_*(h_1-v_1)=h_2-v_2$. The
class $2h_1-3v_1$ is given by the proper transform of the section
of $Y_1$ by $A=0$. Since $f$ maps this section of $Y_1$ into the
section of $Y_2$ by $x_1=0$, it implies $f_*(2h_1-3v_1)=h_2$.
\enddemo

\heading \S3. The adjacency graphs \endheading

\subheading{3.1. Involutions on unimodular lattices} Let $L$ be an
unimodular lattice and $c:L\to L$ its involution (that is an
isometry with $c^2=\id$). Denote by $L_\pm=\Ker (1\mp c)$ the
eigenlattices and by $\frak L_\pm$ their discriminant groups. If
$L$ is odd, then suppose that there exists a characteristic
element $w_L$ belonging either to $L_+$ or to $L_-$ (recall that
{\it characteristic} are the elements $w\in L$ such that $x^2=\la
w,x\ra\mod 2$ for any $x\in L$). This hypothesis allows us to
extend the usual definition of the discriminant quadratic forms
$\frak q_\pm:\frak L_\pm\to\Q/2\Z$ to the case of odd $L$: we put
$\frak q_\pm(x_\pm+L_\pm)=x_\pm^2+\la w_L,x_\pm\ra\mod 2\Z$ for
any $x_\pm+ L_\pm\in\frak L_\pm$, where $w_L\in L_\pm$ is a
characteristic element of $L$ ($w_L\in L_\pm$ is well-defined
modulo $2L_\pm$, thus, $\frak q_\pm$ is independent of its
choice).

As is well-known, $\frak L_\pm$ are $2$-periodic finite groups
with anti-isometric quadratic forms. The group isomorphisms
between $\frak L_\pm$ and $\frak L(c)=L/(L_++L_-)$ providing such
an anti-isometry are induced by the orthogonal projection
$p_\pm\:L\to L_\pm\otimes\Q$. If, for example, $w_L\in L_-$ then
$L_+$ is even and the above isomorphism between $\frak L_+$ and
$\frak L(c)$ converts $\frak q_+$ into $\frak q_c: \frak L(c)\to
\Q/2\Z$ given by $\frak q_c(x+L_++L_-)=\frac12({x^2+\la
x,cx\ra})$, while $\frak q_-$ is converted into $\frac12(x^2-\la
x,cx\ra)+\la x,w_L\ra=-\frak q_c(x+L_++L_-)\mod 2\Z$.

The rank of the (isomorphic) $2$-periodic groups $\frak L_\pm$ and
$\frak L(c)$ will be denoted by $d(c)$.

We distinguish three types of elements in $L_\pm$. An element
$h\in L_\pm$ is called {\it even} if $(h,l)$ is even for any $l\in
L_\pm$, otherwise, $h$ is called {\it odd}. If $h$ is even, then
$\frac12{h}$ defines an element of $\frak L_\pm$ and we call $h$ a
{\it Wu element} if $\frac12{h}\in\frak L_\pm $ is the
characteristic element of the bilinear form $\frak b_\pm:\frak
L_\pm\times\frak L_\pm\to\Q/\Z$. As is known (and easy to check),
$h\in L_\pm$ is even if and only if $h\in (1\pm c)L$. In addition,
an even $h\in L_\pm$ is a Wu element if and only if $x^2+\la
x,cx\ra=\la x,h\ra\mod 2$ for any $x\in L$.

A lattice involution $c \:L\to L$ is called {\it even}, or {\it of
type} I, if $\la x,cx\ra+x^2=0\mod2\Z$ for any $x\in L$.
Otherwise, it is called {\it odd}, or {\it of type} II. Note, that
$c$ is even if and only if the (isomorphic) bilinear forms $\frak
b_\pm:\frak L_\pm\times\frak L_\pm\to \Q/\Z$ are even (that is
$\frak b_\pm(x,x)$ vanishes, or in the other words, the
characteristic element of $\frak b_\pm$ vanishes, or equivalently,
$\frak q_\pm(x)$ is integral for any $x\in\frak L_\pm$).

\subheading{3.2. Adjacency graphs $\Gamma_{d,n}$ and lattice
graphs $\Gamma_L$} For any $d,n\ge1$, one can consider the graph
$\Gamma_{d,n}$ whose vertices are the coarse deformation classes
of real nonsingular hypersurfaces of degree $d$ in $P^{n+1}$ and
the edges are the coarse pseudo-deformation classes in
$\Delta_0(\R)$: according to Lemma 2.1, connected components of
$\Delta_0(\R)$ are topological manifolds; therefore, such a
component and its orbit under the action of $PGL(n+1,\R)$ is
adjacent to at most two coarse deformation classes of nonsingular
hypersurfaces, which are, by definition, the vertices joined by
the edge. As it follows from Corollary 2.8, if $n=2m$, then this
graph has no {\it graph-loops}, that is edges joining a vertex
with itself, and moreover, since connected components of
$\Delta_0(\R)$ are co-oriented and this co-orientation is
preserved by the actions of $PGL(n+1,\R)$, we can orient the graph
pointing the edge from the ascendant side vertex to the opposite
side vertex. For any $d,n\ge1$, the graph $\Gamma_{d,n}$ is
connected as a non-oriented graph.

Replacing the vertices of $\Gamma_{d,n}$ by the middle homology
lattices of the corresponding hypersurfaces and the edges by
Picard-Lefschetz transformations as in the beginning of Section
2.4 we come naturally to the following general notion of the
lattice graph, $\Gamma_{L}$, for an arbitrary lattice $L$.

The automorphism group $\Aut(L)$ acts by conjugation on the set of
involutions $c:L\to L$ and on the set of pairs $(c,x)$, where
$c\:L\to L$ is an involution and $x\in L_-(c)$. We denote by $V_L$
the set of the orbits of the first action, and by $[c]\in V_L$ the
orbit of an involution $c:L\to L$. We let $[c,x]$ denote the orbit
of $(c,x)$ under the second action, and by $E_L$ the set of orbits
of the pairs $(c,x)$ with $x^2=-2, x\in L_-(c)$.

The {\it lattice graph} $\Gamma_L$ is an oriented graph with the
vertex set $V_L$ and the edge set $E_L$, where an edge $[c,v]\in
E_L$ has $[c]$ as its initial vertex and $[c_v]$ with
$c_v=c\,s_v=s_vc$ as its terminal vertex. The elements $[c]\in
V_L$ and $[c_v]\in V_L$ (respectively, the involutions $c$, $c_v$)
are called {\it adjacent}, if they are joined by an edge in
$\Gamma_L$.

The map which associates with a real nonsingular hypersurface $X$
of degree $d$ in $P^{2m+1}$ the lattice $L=(-1)^{m+1}H_{2m}(X)$
and the involution $(-1)^{m+1}c:L\to L$, where $c$ is induced by
the complex conjugation, defines a morphism of oriented graphs
$\Phi_{d,2m}\:\Gamma_{d,2m}\to\Gamma_{L}$. (It may be worth to
notice that $\Gamma_L$ splits as a disjoint union of subgraphs
$\Gamma_{L,p}$, corresponding to a fixed value
$p=\sigma_+(L_+(c))$, and the image of $\Phi_{d,2m}$ is contained
in the subgraph $\Gamma_{L,p}$, where $p=\frac12(\sigma_+(L)-1)$
for odd $m$ and $p=\frac12\sigma_+(L)$ for even.)

An orthogonal pair $v_1,v_2\in L_-(c)$, $v_i^2=-2$, defines two
paths in $\Gamma_L$ which connect the vertices $[c]$ and
$[c\,s_{v_1}s_{v_2}]=[c\,s_{v_2}s_{v_1}]$: one paths consists of
the edge $[c,v_1]$ followed by the edge $[c_{v_1},v_2]$, and the
other one consists of $[c,v_2]$ followed by $[c_{v_2},v_1]$. The
group $\Aut(L)$ acts on the set of such triples $(c,v_1,v_2)$; we
denote by $[c,v_1,v_2]$ the orbit of $(c,v_1,v_2)$, and by $C_{L}$
the set of all such orbits. The graph-cycle in $\Gamma_L$ that we
obtain by following the first path, from $[c]$ along edges
$[c,v_1]$, $[c_{v_1},v_2]$ and then returning along the other path
will be called {\it an elementary cycle}; it depends obviously
only on the orbit $[c,v_1,v_2]$, which will be used as a label for
this cycle. The vertex $[c]$ will be called {\it the origin of the
 elementary cycle $[c,v_1,v_2]$}.

\subheading{3.3. Characterization of the edges in $\Gamma_L$} For
$x,y\in L_-(c)$ let us write $x\sim_c y$  if either \roster \item
$x$ and $y$ are odd, \item $x$ and $y$ are Wu elements (in
particular the both are even), \item $x$ and $y$ are even and both
are not Wu elements.
\endroster

\tm{Proposition 3.1} If $[c,v]$ is an edge of $\Gamma_L$, then the
ranks $d(c)$ and $d(c_v)$ at the endpoints, $[c]$ and $[c_v]$, of
$[c,v]$ are related as follows
 \roster \item if $v$ is odd, then $d(c_v)=d(c)+1$
\item if $v$ is even, then $d(c_v)=d(c)-1$, \item $v$ is a Wu
element with respect to $c$ if and only if $[c_v]$ is even.
\endroster

In particular, if the involution $c_v$ is even, then
$d(c_v)=d(c)-1$.
\endtm

\demo{Proof}
 $L_-(c_v)\oplus \Z v= L_-(c)$
if  $v\in L_-(c)$ is even; otherwise $L_-(c_v)\oplus \Z v$ is a
subgroup of index $2$ in $L_-(c)$. This implies (1) and (2) and
allows to consider $\frak L(c_v)$ as an index $2$ subgroup of
$\frak L(c)$ if $v$ is even.

Assume that $v$ is even. Since $(\la x,c_vx\ra=\la x,cx+\la
cx,v\ra v\ra=\la x,cx\ra+\la cx,v\ra\la x,v\ra=\la x,cx\ra-\la
x,v\ra^2$ has the same parity as $\la x,cx\ra+\la x,v\ra$, the
characteristic classes $\operatorname{wu}(c)\in \frak L(c)$ and
$\operatorname{wu}(c_v)\in \frak L(c_v)\subset \frak L(c)$ of
involutions $c$ and $c_v$ satisfy the relation
$\operatorname{wu}(c_v)=\operatorname{wu}(c)+v$.\qed\enddemo

The graph $\Gamma_L$ has obviously no graph-loops but in general
it may have {\it multiple edges} (two or more edges connecting the
same pair of vertices). In our case of interest, that is for the
K3-lattice $L=3U\oplus 2E_8$, as it follows from Proposition 3.3
below, the graph $\Gamma_L$ does not contain multiple edges.

The proof of Proposition 3.3 and, moreover, the classification of
edges in the case of K3-lattice makes use of the following result
that can be found in \cite{N1}.

\tm{Theorem 3.2} Assume that even non-degenerate lattices $L_i$,
$i=1,2,$ with $2$-periodic discriminants $\frak L_i$ have the same
inertia indices $\sigma_\pm(L_1)=\sigma_\pm(L_2)$ and isomorphic
finite quadratic forms $\frak q_i:\frak L_i\to\Q/2\Z$. Then $L_i$
are isomorphic, and moreover, any isomorphism $\frak L_1\to \frak
L_2$ of quadratic forms $\frak q_i$, $i=1,2$, is induced by an
isomorphism $L_1\to L_2$ if any of the following conditions is
satisfied \roster\item $L_i$ are indefinite lattices; \item $L_i$
are definite lattices of the rank $\le2$. \qed\endroster
\endtm

\rk{Remark} The case (1) in Theorem 3.2 is non-trivial and it is
found in \cite{N1}, whereas (2) is straightforward, because the
only definite lattices of ranks 1 and 2, with $2$-periodic
discriminants are $\la\pm2\ra$ and $2\la\pm2\ra$ respectively.
\endrk

\tm{Proposition 3.3} Assume that $L=3U\oplus 2E_8$ and that the
both eigenlattices $L_\pm=L_\pm(c)$ of a lattice involution
$c:L\to L$ are indefinite. Then for any $v_i\in L_-$, $v_i^2=-2$,
$i=1,2$, the following conditions are equivalent
 \roster\item $[c,v_1]=[c,v_2]$;
\item $[c,v_1]$ and $[c,v_2]$ have the same endpoints,
$[c_{v_1}]=[c_{v_2}]\in V_L$; \item $v_1\sim_c v_2$.
\endroster
\endtm

\demo{Proof} $(1)$ trivially implies $(2)$, and by Proposition
3.1, $(2)$ implies $(3)$.

Assume $(3)$, consider the orthogonal complements of $v_i$,
$L_-^{v_i}=\{x\in L_-\,|\,\la x,v_i\ra=0\}$, and put $\frak
L^i=\discr(L_-^{v_i})$. By definition of $\sim_c$ the elements
$v_1,v_2\in L_-$ are either both even or both odd.

First, restrict our attention to the case of even $v_i$. Under
this assumption, there are orthogonal direct sum decompositions,
$L_-=L_-^{v_i}\oplus \Z v_i$ which give orthogonal direct sum
decompositions of the discriminants, $\frak L_-=\frak
L^i\oplus\discr\la-2\ra$. It implies that $\frak L^i$ have the
same rank and Brown invariant (it is sufficient to apply the
additivity properties of these invariants). They have also the
same parity: indeed, it is sufficient to notice that $\frak
L^i=\frak L(c_{v_i})$ and apply Proposition 3.1. Therefore, $\frak
L^i$ are isomorphic as quadratic spaces (see, for example,
\cite{W1} and \cite{GM}). Furthermore, due to Theorem 3.2 the
lattices $L_-^{v_i}$ are isomorphic. The established isometry
$L_-^{v_1} \to L_-^{v_2}$ extends to an automorphism $\phi_-$ of
$L_-=L_-^{v_i}\oplus\Z v_i$ sending $v_1$ to $v_2$. Applying
Theorem 3.2 to $L_+$ we can extend the induced by $\phi_-$
automorphism of $\frak L_+=-\frak L_-$ to an automorphism $\phi_+$
of $L_+$. The pair $(\phi_+,\phi_-)$ yields an automorphism of
$(L,c)$ sending $v_1$ to $v_2$.

In the case of odd $v_i$, we have $\frak L^i=-\discr(L_+\oplus\Z
v_i)=\frak L_-\oplus\discr\la 2\ra$. The underlying isomorphism
$\frak L^i\to -\discr(L_+\oplus\Z v_i)$ transforms into $v_i$ an
element $w_i\in L_-^{v_i}$ such that $\frac12(v_i+w_i)$ belongs to
$L_-$. Thus, there is an isometry $\frak L^1\to\frak L^2$ which
transforms such a $\frac12 w_1+L_-^{v_1}$ into $\frac12
w_2+L_-^{v_2}$. Therefore, its lift $\phi: L_-^{v_1}\to L_-^{v_2}$
given by Theorem 3.2 extends to an isometry $\phi_1: L_-\to L_-$
sending $v_1$ to $v_2$. Applying once more Theorem 3.2 like in the
previous case we extend $\phi_1$ to an isometry of $L$.
\qed\enddemo

\subheading{3.4. The K3-graph} It is convenient to understand by a
K3-surface any (not necessarily projective and not necessarily
algebraic) nonsingular compact simply-connected complex surface
with the trivial first Chern class and to speak on real
deformations of real K3-surfaces in a sense of Kodaira-Spencer.
Since K3-surfaces have $h^{0,2}=1$ and the complex conjugation
involution maps $H^{0,2}$ to $H^{2,0}$, the involutions induced by
the complex conjugation on the K3-lattice $L=3U\oplus 2E_8$ have
$\sigma_+(L_+)=1$ and $\sigma_+(L_-)=2$, so that these involutions
belong to the component $\Gamma_{L,1}$ of $\Gamma_L$. We call
K3-{\it graph} the graph $\Gamma_{L,1}$ with $L=3U\oplus 2E_8$ and
denote this graph, its vertex, edge and elementary cycle sets by,
respectively, $\Gamma_{K3}$, $V_{K3}$, $E_{K3}$, and $C_{K3}$.
Involutions $c$ with $[c]\in V_{K3}$ will be called {\it real
K3-involutions}.

 \midinsert \epsfbox{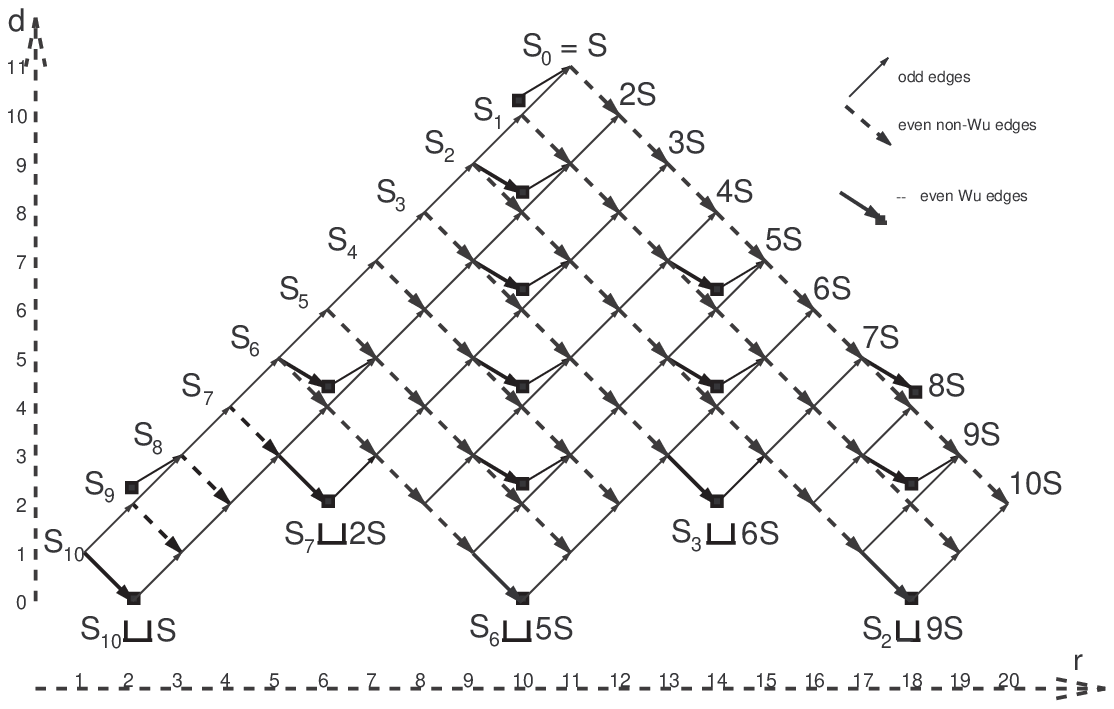}
\topcaption{Figure 1. The K3-graph $\Gamma_{K3}$} \endcaption
\noindent {\eightrm  Vertices denoted by $\blacksquare$ represent
K3-surfaces of type I. Other vertices represent K3-surfaces of
type II. The coordinates are $r=10+\frac12\chi(X(\R))$ and
$d=12-\frac12 \dim H_*(X(\R);\Z/2\Z)$ for $X(\R)\ne\varnothing$.
The case $(r,d)=(10,10)$ corresponds to $X(\R)=\varnothing$ (type
I) and $X(\R)=S_1$ (type II). The case $(r,d)=(10,9)$ corresponds
to $X(\R)= 2S_1$ (type I) and $X(\R)=S_2\+ S$ (type II). Other
surfaces have $X(\R)=S_p\+ qS$, where $p\ge 0$ and $q\ge 0$ are
uniquely determined by $(r,d)$.}
\endinsert

\tm{Theorem 3.4} \roster \item A vertex $[c]\in V_{K3}$ is
determined by the isomorphism type of $L_+(c)$. \item The
isomorphism type of $L_+(c)$ for $[c]\in V_{K3}$ is determined by
the rank of $L_+(c)$ and the rank and parity of its discriminant.
\item The graph $\Gamma_{K3}$ is like is shown on Figure 1.
\endroster
\endtm

\demo{Proof} Assertions (1) and (2), as well as a list of
necessary and sufficient conditions on the rank of $L_+(c)$ and
the rank and parity of its discriminant, are contained in
\cite{N1}.  This gives the set of vertices of $\Gamma_{K3}$. All
the restrictions on the edges follow from Propositions 3.1 and
3.3, and their existence follows from Propositions 4.2--4.3
applied to the list of the real K3-involutions given in Tables
2--3 of the Appendix. (Note that the existence of all but one edge
can be deduced also from \cite{Itenberg}, which contains a
classification of adjacency-edges in the case of $2$-polarized
K3-surfaces.) \qed
\enddemo

One can observe furthermore that the cycles $[c,v_1,v_2]\in
C_{K3}$, where $v_1$ is an odd and $v_2$ is an even nodal class,
$v_i\in L_-(c)$, form a basis for $1$-homology of $\Gamma_{K3}$.
We will call them {\it basic cycles}.

It may be worth to mention (although not really essential for us)
that the topological type of the real locus, $X(\R)$, of a
K3-surface $X$ determines the rank $r$ of $L_+$ and the rank $d$
of its discriminant (see Figure 1; more details and the references
can be found, for example, in \cite{DIK1}). In a number of cases
the topological type of $X(\R)$ determines as well the parity of
$c$, then it defines a certain vertex $[c]\in V_{K3}$, and,
moreover, the topological type is uniquely defined by such a
vertex.
 This is the first kind of vertices. To identify a vertex
$[c]\in V_{K3}$ of the second kind we need to determine the parity
of the involution. Therefore, when it is appropriate, we will use
notation $[T]$, for the vertices of the first kind and for the odd
vertices of the second kind. For the even vertices of the second
kind we use notation $[T]_I$. Here $T$ is a topological type from
the list of real K3-surfaces (see Figure 1), which can be
$S_p\dsum qS$ (disjoint union of an orientable surface of genus
$p$ and $q$ spheres), or $2S_1$ (disjoint pair of tori), or
$\emptyset$.

\subheading{3.5. The K4-graph} The graph $\Gamma_{3,4}$ generated
by cubic fourfolds is our principal object of interest. We call it
{\it K4-graph} and denote this graph, its vertex set, and its edge
set by, respectively, $\Gamma_{K4}$, $V_{K4}$, and $E_{K4}$.

Proposition 2.3 yields a bijection between $E_{K4}$ and the set of
coarse deformation classes of bi-degree $(2,3)$ real K3-surfaces
in $P^4$. On the other hand, as is shown in \cite{N1} the latter
set is in one-to-one correspondence with the set of isomorphism
classes $[c,h]$ such that $[c]\in V_{K3}$, $h\in L_-(c)$, and
$h^2=6$. We will use $[c,h]$ as a label for the corresponding
K4-edge and denote by $w[c,h]$ and $w_+[c,h]$ its initial and
terminal vertices.

The classification of isomorphism classes $[c,h]$ such that
$[c]\in V_{K3}$, $h\in L_-(c)$, and $h^2=6$ which is given in
\cite{N06} can be reformulated as the following classification of
K4-edges.

\tm{Theorem 3.5} Edges $[c,h_1], [c,h_2]\in E_{K4}$ coincide if
and only if $h_1\sim_c h_2$. In particular, for any K3-vertex
$[c]$, there exist at most three K4-edges $[c,h]$, namely, with
$h\in L_-$ being (1) an odd element, (2) a Wu element, or (3) even
but not Wu element.\qed\endtm

These three edges are characterized by the following relation
between the discriminant ranks $d(w)$ and $d(w_+)$ of the
K4-lattice involutions $c_{M(w)} : M(w)\to M(w)$ and $
c_{M(w_+)}\:M(w_+) \to M(w_+)$ of the cubic fourfolds
corresponding to the adjacent pair of vertices $w=w[c,h]$ and
$w_+=w_+[c,h]$.

\tm{Corollary 3.6} For any K4-edge $[c,h]$ \roster \item if $h$ is
odd, then $d(w_+)=d(w)+1$, \item if $h$ is even, then
$d(w_+)=d(w)-1$, \item the involution in the lattice $M(w_+)$ is
even if and only if $h$ is a Wu element with respect to $c$.
\endroster
\endtm

\demo{Proof} According to Corollary 2.10, $M(w)$ is canonically
identified as a group with $L\oplus \Z$, and under this
identification the involution $c_{M(w)}$ coincides with
$-c\oplus\id$. Furthermore, by Lemma 2.12 the vanishing class
$v\in M_+(w)$ is of the same parity as $h$, namely $v=h+2e$ where
$e$ is a generator of $\Z$. To prove assertions (1) and (2) it
remains to notice that the involution $c_{M(w_+)}$ in $M(w_+)$ is
equal to $s_v\circ c_{M(w)}$, to use Lemma 2.6 and to apply
Proposition 3.1. To prove assertion (3) one needs in addition to
check that $c_{M(w_+)}$ is even if and only if  $h$ is a Wu
element of $c$. The latter equivalence follows from a $v$-twist
relation $\la x+ne,s_v c_{M(w_+)}(x+ne)\ra+\la x+ne,s_v(x+ne)\ra
=x^2+\la x,cx\ra+\la x,h\ra\mod 2\Z$ for any $x\in L$ and
$n\in\Z$, which is straightforward. \qed
\enddemo

\subheading{3.6. K4-cycles} The automorphism group $Aut(L)$ of the
K3-lattice $L$ naturally acts on the set of triples $(c,h,v)$,
where $c\in V_{K3}$, $h,v\in L_-(c)$, $h^2=6$, $v^2=-2$, and
$v\bot h$. Let $[c,h,v]$ denote the orbit of a triple $(c,h,v)$
and $C_{K4}$ the set of all such orbits.

The following assertion is well known. We give its proof, since
did not find a straightforward reference.

\tm{Theorem 3.7} Any $[c,h_1,v_1]\in C_{K4}$ can be represented by
the complex conjugation, hyperplane section, and exceptional
divisor of the nonsingular model of a real bi-degree $(2,3)$
complete intersection K3-surface in $P^4$ with a single node.
\endtm

\demo{Proof} Pick a generic $\omega=\omega_++i\omega_-$ with
$\omega_+\in L_+\otimes\R$, $\omega_-\in L_-\otimes\R$,
$\omega_+^2=\omega_-^2>0$, and  $(\omega, h)=(\omega,v)=0$
(generic in a sense that there is no $l\in L$ orthogonal to
$\omega$ other than linear combinations of $h$ and $v$). Due to
the surjectivity of the period map, there exists a marked
K3-surface $Y$ equipped with an anti-holomorphic involution $\tau
:Y\to Y$ and an isometry $\phi:H^2(Y;\Z)\to L$ such that
$c\circ\phi=\phi\circ\tau_*$ and
$\phi^{-1}\omega\otimes\C=H^{2,0}(Y)$ (to take into account the
anti-holomorphic involution, one can apply, for example,
\cite{DIK1}, Theorem 13.4.3). Due to the prescribed isomorphism
type of the lattice generated by $h$ and $v$ and the generic
choice of $\omega$, the Picard group $\Pic Y\subset H^2(Y;\Z)$ is
generated by the fundamental classes $v_1,v_2$ of two
$(-2)$-curves with $(v_1, v_2)=4$. Therefore, one can adjust the
marking $\phi$ to have $\phi(v_2)=v$ and $\phi(v_1+2v_2)=h$. Then,
$H=v_1+2v_2$ is a big nef divisor, and since, in addition,
$E^2\ne 0$ for any $E\in\Pic(Y)$, it follows from \cite{StD},
Propositions 2.6 and Theorem 5.2, that the linear system $\vert
H\vert $ has no fixed components or fixed points and, moreover,
defines a degree one map onto a normal surface with simple
singularities (in our case, the only singularity is the node
representing $v_2$).

It remains to notice that, according, for example, to \cite{StD}
Theorem 6.1, any degree $6$ normal K3-surface in $P^4$ is a
complete intersection of a quadric with a cubic hypersurface.
\qed\enddemo

Recall that the central projection correspondence associates to a
uni-nodal K3-surface $Y_1$, which is a degree $(2,3)$ complete
intersection in $P^4$, a bi-nodal cubic fourfold $X_0$ with one of
the nodes, $x_1\in X_0$, being marked. Marking the other node,
$x_2\in X_0$, yields another uni-nodal K3-surface, $Y_2$, related
similarly to $[c,h_2,v_2]$, where $h_2=2h_1-3v_1$, $v_2=h_1-2v_1$
(see Proposition 2.13). The involution $C_{K4}\to C_{K4}$,
$[c,h,v]\mapsto [c,2h-3v,h-2v]$, will be called {\it the flip
involution}.

In the space $P^N$ of all cubic fourfolds, the bi-nodal cubics
represent the points of transversal self-intersection of the
discriminant hypersurface $\Delta\subset P^{N}$ (cf. 2.1) and the
two (intersecting transversally) branches of $\Delta$ correspond
to the two nodes. The branches are real if the nodes are real. If
the nodes are real, then each of these real branches is divided by
their intersection into a pair of ``halves'', which are the
adjacent pseudo-deformation classes representing K4-edges
$[c,h_i]$ and $[c_{v_i},h_i]$ (see Figure 2).
 These edges form a cycle connecting the four vertices of
 $\Gamma_{K4}$ presented by the coarse deformation components in
 $P^{N}(\R)$
which are locally separated by the pair of branches of
$\Delta(\R)$. Being oriented in accordance with the orientation of
the edge $[c,h_i]$, such cycle will be labelled by $[c,h_i,v_i]$.
Note that the flip $[c,h_1,v_1]\mapsto [c,h_2,v_2]$ just changes
the orientation of the cycle, like the involution $C_{K3}\to
C_{K3}$, $[c,v_1,v_2]\mapsto [c,v_2,v_1]$ changes the orientation
of the elementary K3-cycles.

\midinsert \epsfbox{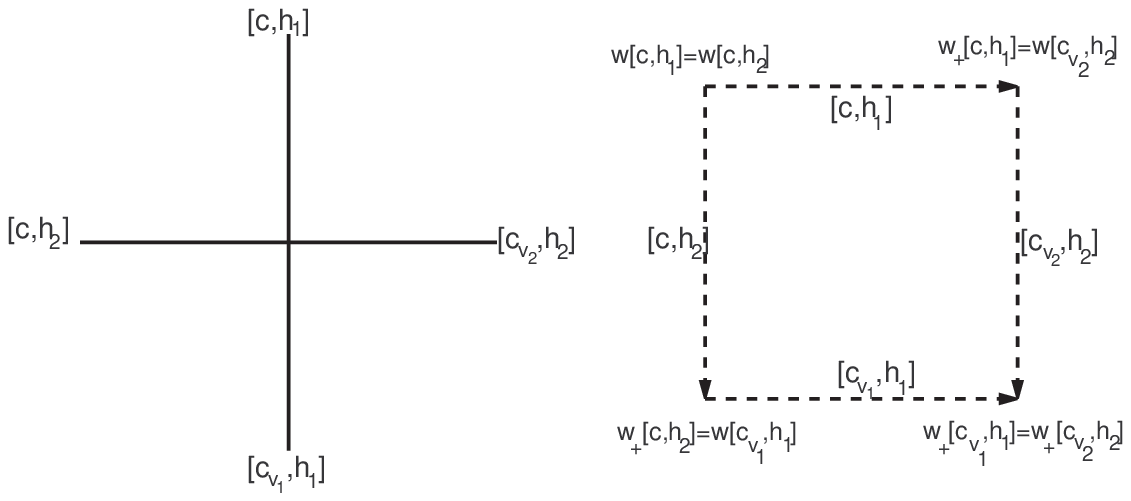} \botcaption{Figure 2} {\eightrm
Self-intersection of the discriminant $\Delta(\R)$ and the
corresponding K4-cycle}\endcaption\endinsert

Summarizing, we obtain the following result.

\tm{Proposition 3.8} For any real K3-lattice involution $c\:L\to
L$ and classes $h_1,v_1\in L_-(c)$, $h_1^2=6$, $v_1^2=-2$,
$h_1\bot v_1$, there exists a cycle in $\Gamma_{K4}$ formed by
four edges $[c,h_1]$, $[c,h_2]$, $[c_{v_1},h_1]$, and
$[c_{v_2},h_2]$, where $h_2=2h_1-3v_1$, $v_2=h_1-2v_1$.

In the other words, the following endpoints of these edges
coincide
$$\align
w[c,h_1]&=w[c,h_2]\\
w_+[c,h_1]&=w[c_{v_2},h_2]\\
w_+[c,h_2]&=w[c_{v_1},h_1]\\
w_+[c_{v_1},h_1]&=w_+[c_{v_2},h_2] \qed\endalign$$
\endtm

\heading \S4. The graph isomorphism
$\Gamma_{K4}^*\cong\Gamma_{K3}^*$
\endheading

\subheading{4.1. The graph morphism $\Gamma_{K4}^*\to
\Gamma_{K3}^*$}
 We say that vertices $[c]\in V_{K3}$ and $w\in V_{K4}$
correspond to each other if there is an anti-isometry between
their eigenlattices, $M_-(w) =-L_+(c)$. We say also that there is
correspondence between edges $[c,v]\in E_{K3}$, $[c,h]\in E_{K4}$
if $h\sim_c v$. Cycles $[c,v_1,v_2]\in C_{K3}$ and $[c,h,v]\in
C_{K4}$ correspond to each other if $h\sim_cv_1$ and $v\sim_c
v_2$. Note that if cycles $[c,v_1,v_2]$ and $[c,h,v]$ correspond
to each other, then after their orientation is reversed they also
correspond to each other, i.e., $[c,v_2,v_1]$ corresponds to
$[c,2h-3v,h-2v]$. This is because $2h-3v\sim_c v\sim_c v_2$ and
$h-2v\sim_c h\sim_c v_1$.

A vertex, an edge, or a cycle (in $\Gamma_{K3}$ or $\Gamma_{K4}$)
will be called {\it regular} if it has a corresponding one (in the
other graph), otherwise it will be called {\it irregular}. The
sets of regular vertices, edges, and cycles in $\Gamma_{K3}$ (or
$\Gamma_{K4}$) will be denoted $V_{K3}^*$, $E_{K3}^*$, and
$C_{K3}^*$ (respectively, $V_{K4}^*$, $E_{K4}^*$, and $C_{K4}^*$).

Theorem 3.4(1-2) implies the uniqueness of K3-vertex $[c]$ that
corresponds to a vertex $w\in V_{K4}^*$ and Proposition 3.3
implies similarly the uniqueness of $[c,v]\in V_{K3}^*$ that
corresponds to $[c,h]\in V_{K4}^*$. Let $F_V\:V_{K4}^*\to
V_{K3}^*$ and $F_E\:E_{K4}^*\to E_{K3}^*$ denote the maps defined
by this correspondence.

\tm{Proposition 4.1} $\Gamma_{K3}^*=(V_{K3}^*,E_{K3}^*)$ and
$\Gamma_{K4}^*=(V_{K4}^*,E_{K4}^*)$ are subgraphs of $\Gamma_{K3}$
and $\Gamma_{K4}$ spanned by the vertex sets $V_{K3}^*$ and
$V_{K4}^*$, and the pair of maps $F=(F_V,F_E)$ is a morphism of
oriented graphs, $F\:\Gamma_{K4}^*\to\Gamma_{K3}^*$. More
precisely, \roster \item the initial vertex $w[c,h]$ is regular
for any edge $[c,h]\in E_{K4}$; the terminal vertex $w_+[c,h]$ is
regular if and only if the edge $[c,h]$ is regular; \item an edge
$[c,v]\in E_{K3}$ is regular if and only if its both endpoints are
regular; \item $F$ preserves incidences of edges and vertices and
the order of vertices of an edge; that is, if $[c,h]\in E_{K4}^*$
and $F_E([c,h])=[c,v]$, then $F_V(w[c,h])=[c]$ and
$F_V(w_+[c,h])=[c_v]$; \item $F_E$ is bijective; \item $F$ sends a
regular cycle $[c,h,v]\in C_{K4}^*$ into the corresponding one,
$[c,v_1,v_2]\in C_{K3}^*$, $h\sim_c v_1$, $v\sim_c v_2$ (in
particular, all vertices and edges of the corresponding to each
other cycles are regular).
\endroster
\endtm

\demo{Proof} By Proposition 2.11, $L_+(c) =-M_-(w)$ for any
$[c,h]\in E_{K4}$, so $w=w[c,h]$ is regular and $F(w)=[c]$.

Assume that $[c,h]\in E_{K4}$ is a regular edge and consider $v\in
L_-(c)$ such that $v^2=-2$ and $h\sim_c v$. If $h,v\in L_-(c)$ are
odd, then $L_+(c_v)=L_+(c)\oplus\la-2\ra$ and
$M_-(w_+)=M_-(w)\oplus\la2\ra$, which implies $F(w_+)=[c_v]$. If
$h,v\in L_-(c)$ are even, then a similar decomposition for the
complementary lattices, $L_-(c)=L_-(c_v)\oplus\la-2\ra$ and
$M_+(w)=M_+(w_+)\oplus\la2\ra$, implies that the discriminants of
$M_\pm(w_+)$ and $-L_\mp(c_v)$ have the same rank and the same
Brown invariant. By Proposition 3.1 and Corollary 3.6, their
discriminant forms have the same parity. Thus, Theorem 3.2 applied
to the even lattices $M_-(w_+)$ and $-L_+(c_v)$ shows that they
are isometric.

Conversely, if $w_+[c,h]$ is regular and corresponds to $[c']$,
then $L_+(c')=-M_-(w_+)$, where the latter, on the other hand, is
either isomorphic to the lattice $-(M_-(w)\oplus\la2\ra)$, or
contains it as a subgroup of index 2. Therefore, $c=s_v\circ c'$,
where $v\in L_+(c')$ corresponds to the second summand in
$-(M_-(w)\oplus\la2\ra)\subset L_+(c')$. Moreover, according to
Proposition 3.1 and Corollary 3.6, $h\sim_c v$. This completes the
proof  of (1) and (3). Furthermore, (2) is their straightforward
corollary.

Surjectivity of $F_E$ holds by its definition and injectivity is
due to Theorem 3.5, which yields (4).

For (5), note that $[c,h]$ corresponds to $[c,v_1]$ because
$h\sim_c v_1$, $[c,2h-3v]$ corresponds to $[c,v_2]$ because
$2h-3v\sim_cv$, $[c_{v},h]$ corresponds to $[c_{v_2},v_1]$ because
$[c_v]=[c_{v_2}]$ by Proposition 3.3, and similarly,
$[c_{h-2v},2h-3v]$ corresponds to $[c_{v_1},v_2]$. \qed\enddemo

\subheading{4.2. Regularity of K3 and K4 edges} A simultaneous
description of K3 and K4 edges in the following two lemmas implies
their regularity with two exceptions: one K3-edge and one K4-edge.

We refer to the Appendix for the notation of the standard lattices
which appear in the direct sum decomposition of $L_\pm(c)$. By the
{\it diagonal component} of $L_-(c)$ for $[c]\in V_{K3}$ we mean
the component $s\,\la2\ra\oplus t\,\la-2\ra$ in the decomposition
of $L_-(c)$, as presented in Tables 2--3 in the Appendix.  Note
that $c$ is an even involution if and only if $t=s=0$. In
particular, the diagonal component vanishes for all the
involutions in Table 3. An opposite kind of extremal case is
$L_-(c)=s\,\la2\ra\oplus t\,\la-2\ra$, which corresponds to
$[c]=[kS]$, $1\le k\le10$, $s=2$ and $t=10-k$, as is seen from
Table 2. It is also obvious from this table that K3-vertices
$[10S]$ and $[8S]_I$ are terminal, i.e., do not have outgoing
edges, since $L_-([10S])=2\,\la2\ra$ is positive, and all $x\in
L_-([8S]_I)=U(2)$, have $x^2$ divisible by $4$.

\tm{Proposition 4.2} For any given $[c]\in V_{K3}$ and
$n\in\{0,1\}$, the following conditions are equivalent:
\roster\item there exists an odd element $x\in L_-(c)$,
$x^2=8n-2$; \item $[c]$ differs from $[kS]$, $k\ge1$, and
$[8S]_I$. \endroster In particular, all the odd edges, $[c,v]\in
E_{K3}$ and $[c,h]\in E_{K4}$, are regular.
\endtm

\demo{Proof} All elements of
$L_-([kS])=2\,\la2\ra\oplus(10-k)\la-2\ra$ and $L_-([8S]_I)=U(2)$
are even, so (1) implies (2). Assuming (2), we can observe in
Tables 2--3 that $L_-(c)$ contains as a direct summand one of the
following three lattices: $U$, or $\la2\ra\oplus E_8$, or
$U(2)\oplus D_4$. In the first case, $x=(1,k)\in U$ is an odd
element for any $k\in\Z$, and in particular for $k=4n-1$, which
gives $x^2=2k=8n-2$.
 In the second case, $x=x_1\oplus
x_2\in\la2\ra\oplus E_8$ is an odd element, if $x_1=2n\,e_+$,
$x_2=(2n-1)e_-$, where $e_\pm^2=\pm2$, and $n\in\Z$. It gives
$x^2=8n-2$.
 In the third case, element $x=x_1\oplus x_2\in
U(2)\oplus D_4$ is odd with $x^2=4k-2$ if $x_1=(1,k)$ in $U(2)$
and $x_2\in D_4$, $x_2^2=-2$. \qed\enddemo

\tm{Proposition 4.3} Let $s\,\la2\ra\oplus t\,\la-2\ra$ be the
diagonal component of $L_-(c)$ for $[c]\in V_{K3}$, and
$n\in\{0,1\}$. Then
 \roster \item there exists an even element $x\in L_-(c)$, with
$x^2=8n-2$ if and only if $t\ge1$;
 \item in the
case of $[c]\ne[kS]$, existence of a Wu element $x\in L_-(c)$,
with $x^2=8n-2$ is equivalent to $s-t=-1\mod4$;
 \item in the case of
$[c]=[kS]$, existence of a Wu element $x\in L_-(c)$, with
$x^2=8n-2$ is equivalent to $s-t=4n-1\mod8$;
 \item
there exists a non-Wu even element $x\in L_-(c)$ with $x^2=8n-2$
if and only if one of the following two conditions is satisfied:
(a) $t>1$, or (b) $t=1$ and $s-t\ne-1\mod4$.
\endroster
\endtm

\demo{Proof} Proving (1)--(3), we consider the decomposition of
elements $x\in L_-(c)$ as $x=x_1\oplus x_2$, where $x_1\in
s\la2\ra\oplus t\la-2\ra$ and $x_2$ is from the complementary
non-diagonal component of $L_-(c)$. Note that $x$ is even in
$L_-(c)$ if and only if $x_2$ is even in its non-diagonal
component (because any element
 $x_1\in s\la2\ra\oplus t\la-2\ra$ is even).
It is straightforward to observe also that $x_2^2$ is divisible by
$4$ if $x_2$ is even (just by checking such a divisibility for the
direct summands $U$, $E_8$, $D_4$, $U(2)$ and $E_8(2)$ that appear
in the non-diagonal component of $L_-(c)$). Thus, if $s=t=0$ then
$x^2=x_2^2\ne 8n-2$. If $s>t=0$, then from Table 2 we can see that
the non-diagonal component of $L_-(c)$ contains only summands $U$
and $E_8$. Observing that even elements of unimodular lattices are
obviously divisible by two, we conclude that $x_2^2$ is divisible
by $8$. Thus, $x^2=x_1^2\mod8$, which cannot be $-2\mod8$ for
$s\le2$, because $x_1^2$ is just the doubled square of an integer
if $s=1$ and the doubled sum of squares of two integers if $s=2$.
 This proves that the condition in (1) is necessary.
To show its sufficiency, we observe from Table 2 that if $t\ge1$,
then either $s\ge1$, or $L_-(c)$ contains a summand $U$ orthogonal
to $\la-2\ra$. In the first case,
$x=(k+1,k)\in\la2\ra\oplus\la-2\ra$, is an even element,
$x^2=(4k+2)=8n-2$, if we let $k=2n-1$. If $L_-(c)$ contains
$\la-2\ra\oplus U$, then the sum of a generator of $\la-2\ra$ with
$y=(2,2n)\in U$ is an even element of square $8n-2$. This
completes the proof of $(1)$.

 Since a direct sum decomposition of a
lattice yields a direct sum decomposition of the discriminant
form, $x$ is a Wu element of $L_-(c)$ if and only if $x_1$ and
$x_2$ are Wu elements in the corresponding  components of
$L_-(c)$. It is straightforward to observe that if $x_1$ is a Wu
element, then $x_1^2=2(s-t)\mod 16$, and if $x_2$ is a Wu element,
then $x_2^2$ is divisible by $8$ (indeed, the discriminant form of
$s\la2\ra\oplus t\la-2\ra$ is $s\la\frac12\ra\oplus t\la-\frac12
\ra$, the discriminant forms of $D_4$, $U(2)$, $E_8(2)$ are even,
and those of $U$ and $E_8$ are trivial). This shows that the
conditions $s-t=-1\mod4$ in (2) and $s-t=4n-1$ in (3) are
necessary for existence of a Wu element $x$ with $x^2=8n-2$.

If $s-t=4n-1$ in (3), then for the existence part of the
statement, we need only to observe that
$L_-([7S])=2\,\la2\ra\oplus3\,\la-2\ra$ contains a Wu class $v$,
$v^2=-2$, and that $L_-([3S])=2\,\la2\ra\oplus7\,\la-2\ra$
contains a Wu class $h$, $h^2=6$.

If $s-t=-1\mod 4$ in (2), then we note from Tables 2--3 that
either $L_-(c)$ contains $U$, or it contains $s\la2\ra\oplus
t\la-2\ra\oplus E_8$ with $s=2$. If there is a summand $U$, then
letting $x_2=(2,2k)\subset U$ we obtain
$x^2=x_1^2+x_2^2=2(s-t)+8k$, where $x_1=(1,1,\dots,1)\in
s\la2\ra\oplus t\la-2\ra$. Choosing a suitable $k$, we obtain a Wu
element $x$ with $x^2=8n-2$, if $(s-t)=-1\mod4$. Similarly, we
obtain a Wu element $x$ with $x^2=8n-2$ in the second case. For
example, one can choose $x_1=(3,1,\dots,1)\in s\la2\ra\oplus
t\la-2\ra$, so that $x_1^2=2(s-t)+16$ is either $6$ or $14$ (since
$1\le s\le2$ and $0\le t\le 9$), and then let $x_2=2(\kappa_1
e_1+\kappa_2 e_2)$, for any $e_i\in E_8$, $e_i^2=-2$, $e_1\bot
e_2$ and appropriate $\kappa_i\in\{0,1\}$ to obtain $x^2=8n-2$.

The conditions (a) or (b) in (4) are necessary because in the case
of $t=1$, $s=0$ the discriminant form is either $\la-\frac12\ra$,
or $\la-\frac12\ra\oplus\discr(D_4)$, which implies that non-Wu
even elements $x\in L_-(c)$ have $x^2$ divisible by $4$, in
contradiction to $x^2=8n-2$. Sufficiency of (b) follows from
(1)--(3). In case of (a), the construction of $x$ in (1) gives a
non-Wu even element, because it involves only one $\la-2\ra$
summand. \qed\enddemo

\tm{Corollary 4.4} There exists only one irregular K3-edge,
namely, $[c,v]$ with $[c]=[7S]$, and a Wu element $v\in
L_-(c)=2\,\la2\ra\oplus3\,\la-2\ra$, $v^2=-2$, and only one
irregular K4-edge, namely, $[c',h]$, with $[c']=[3S]$ and a Wu
element $h\in L_-(c')=2\la2\ra\oplus7\la-2\ra$, $h^2=6$.

The terminal vertices of these edges, $[c_v]=[8S]_I\in V_{K3}$ and
$w_+[c',h]\in V_{K4}$ are the only irregular vertices in $V_{K3}$
and $V_{K4}$.
\endtm

\demo{Proof} The criteria of existence for $v$ and $h$ (odd, Wu,
and even non-Wu) are all the same except the one in Proposition
4.3(3). The conclusion about the irregular vertices follows from
Proposition 4.1(1)-(2). \qed\enddemo

\subheading{4.3. Graph $\Gamma_{K4}$ and proof of Theorem 1.1}

\tm{Proposition 4.5} The basic K3-cycles are all regular. More
precisely, for any elementary cycle $[c,v_1,v_2]\in C_{K3}$ such
that $v_1$ is odd and $v_2$ is even there exists an odd element
$h\in L_-(c)$ with $h^2=6$ and $h\bot v_2$, so that $[c,h,v_2]\in
C_{K4}$ corresponds to $[c,v_1,v_2]$.
\endtm

\demo{Proof} Note that $v_1$ is an odd element of $L_-(c_{v_2})$
because of the direct sum decomposition
$L_-(c)=L_-(c_{v_2})\oplus\la-2\ra$, where $\la-2\ra$ is generated
by $v_2$. Proposition 4.2 applied to the involution $c_{v_2}$
implies the existence of an odd element $h\in L_-(c_{v_2})$,
$h^2=6$, and, thus, $[c,h,v_2]\in C_{K4}$ is such as required.
\qed\enddemo

\tm{Proposition 4.6} The mapping $F_V\: V_{K4}^*\to V_{K3}^*$ is
bijective, and thus $F\:\Gamma_{K4}^*\to\Gamma_{K3}^*$ is an
isomorphism of graphs.
\endtm

\demo{Proof} The both graphs, $\Gamma_{K4}^*$ and $\Gamma_{K3}^*$,
are connected, since  they are obtained from connected graphs
$\Gamma_{K4}$ and $\Gamma_{K3}$ by removing irregular vertices and
edges, that is, by Corollary 4.4, just one vertex of valency one
with the adjacent edge from each graph. On the other hand, $F_V$
is surjective, by its definition, $F_E$ is bijective by
Proposition 4.1(4), and $F$ induces an epimorphism at the level of
the first homology by Proposition 4.5. All this together implies
that $F$ is an isomorphism of graphs. \qed\enddemo

\midinsert \topcaption{Figure 3. The K4-graph $\Gamma_{K4}$}
\endcaption\epsfbox{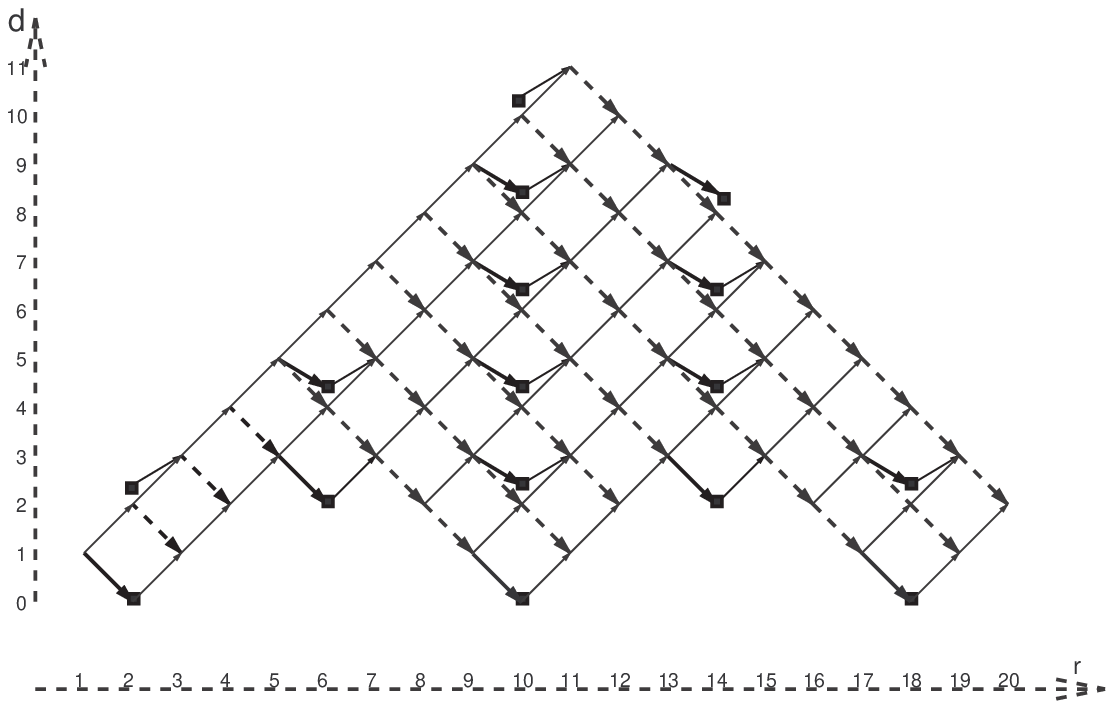}
\eightpoint \noindent The two types of vertices (I and II) and the
three types of edges (odd, Wu, and even non-Wu) on Figure 3 are
marked in the same way as on Figure 1.
\endinsert

\tm{Corollary 4.7} The adjacency graph $\Gamma_{K4}$ for real
nonsingular cubic hypersurfaces in $P^5$ is the graph shown on
Figure 3. In particular, the number of coarse deformation classes
for such hypersurfaces is $75$.
\endtm

\demo{Proof} By Proposition 4.6 and Corollary 4.4, to obtain
$\Gamma_{K4}$ one needs just to remove from $\Gamma_{K3}$ the
irregular K3-edge together with its endpoint $[8S]_I$ and to add
instead an irregular K4-edge with the origin at the K4-vertex
corresponding to $[3S]$. \qed\enddemo

Finally, we will prove Theorem 1.1, which can be reformulated as
follows.

\tm{Theorem 4.8} If the eigenlattices $M_-=\Ker (1+c)$ of two real
nonsingular cubic hypersurfaces in $P^5$ are isomorphic, then
these two hypersurfaces are coarse deformation equivalent. In
particular, the coarse deformation classes of real nonsingular
cubic hypersurfaces in $P^5$ are distinguished by the homological
type of these hypersurfaces.
\endtm

\demo{Proof} Injectivity of $F_V$ implies that regular K4-vertices
$w$ can be distinguished by their eigenlattices
$M_-(w)=-L_+(F_V(w))$, due to Theorem 3.4(1). The irregular
vertex, $w\in V_{K4}\smallsetminus V_{K4}^*$, can be also
distinguished, because its lattice $M_-(w)$ is not isomorphic to
$-L_+(c)$ for all real K3-involutions $c$.

Namely, $-M_-(w)=U(2)+3D_4$, since due to Theorem 3.2 it is the
only even lattice having the required signature, $(1,13)$, and an
even discriminant form of the required rank $d=8$. These values
follows from Corollary 3.6 applied to the irregular K4-edge.
\qed\enddemo

\subheading{4.4. A few related remarks}

\rk{Deformations of cubics}\footnote{We are grateful to
J.Koll\'{a}r for pointing out to us Fujita's results we are using
in this remark.} In the case of nonsingular cubic hypersurfaces of
dimension $\ge 3$ the coarse deformation equivalence can be
reformulated as a deformation equivalence of polarized varieties.

Namely, nonsingular cubic hypersurfaces $X\subset P^{n}$ come with
a natural, projective, polarization $L$, which is of degree $3$
and of $\Delta$-genus $\Delta(X,L)=n-1+L^{n-1}-h^0(X,L)=1$. And
vice versa, according to T.~Fujita \cite{F}
 any polarized
manifold of degree $3$ and $\Delta=1$ is a cubic hypersurface. In
addition, any antiholomoprhic involution $X\to X$ lifts to an
involutive anti-isomorphism of $L$ if $n\ge 4$ (recall that
$X(\R)\ne\varnothing $ for any cubic and $\operatorname{Pic} X=\Z$
for nonsingular cubics in $P^n$ with $n\ge 4$), and therefore any
real structure $c:X\to X$ is induced by the complex conjugation
under a suitable realization of $X$ as a projective cubic if $n\ge
4$.

On the other hand, for nonsingular cubic hypersurfaces it holds
$K+(n-1)L=L$, which implies that $H^i(X;L)=0$ for any $i>0$. By
semi-continuity, it implies the invariance of $\Delta$-genus under
deformations preserving the polarization. As a result, coarse
deformation classes of real nonsingular cubic hypersurfaces of
dimension $n-1\ge 3$ coincide with equivariant deformation classes
of real polarized manifolds with polarization of degree $3$ and
$\Delta$-genus $1$.
\endrk

\rk{Some open problems} When the dimension of the projective space
$P^n$ is even, the group $PGL(n+1,\R)$ is connected, and hence in
this case there is no difference  between coarse deformation and
deformation equivalences as they are defined in Introduction. If
this dimension is odd, the group has two connected components, and
therefore there may be a certain difference between the two
classifications. For the moment, we do not know does it really
happen or not in the case of cubics in $P^5$.

Let us mention a few other open problems concerning real
nonsingular cubics in $P^5$. First, it would be interesting to
give an explicit description of the topology of the real part
$X_\R$  of the cubics in each deformation class.
 It is not difficult to observe, for instance, that K3-vertex
 $[\oo]$ corresponds to the cubics with $X_{\R}=\Rp4\+S^4$, and
$[S]$ corresponds to $X_{\R}=\Rp4$. The simplest candidate to be
proposed in the case of $[S_p\+qS]$ is $X_{\R}=\Rp4\#p(S^2\times
S^2)\#q(S^1\times S^3)$.

The second question is to relate the topology of $X_\R$ with the
topology of the real part of the corresponding Fano varieties,
which are, as is known, are nonsingular for any nonsingular cubic
and hence also depend only on the coarse deformation class of the
cubic. (The case of Fano surfaces of three-dimensional cubics is
worked out in \cite{Kra2}.)
\endrk

\rk{Some other relations between cubic fourfolds and K3-surfaces}
In this paper we related nodal cubic hypersurfaces of dimension
$4$ with K3-surfaces via a central projection from a double point.
As is known, there exist many other connections between them. For
example, if (instead of a double point) a cubic hypersurfaces of
dimension $4$ contains a plane, then the critical locus of the
projection from this plane to a complementary plane is a curve of
degree $6$, so that taking the double covering of the plane
ramified in such a curve one obtains once more a K3-surface.
Another option is to consider the Fano variety of the cubic. In
dimension $4$, this variety (which is always nonsingular) is
deformation equivalent to the symmetric square of a K3-surface.
(It may be worth to mention that these deformation equivalent
varieties are hyperk\"{a}hler and the Fano varieties in question
provide a hypersurface in the corresponding moduli space.)

In fact, the main reason of relating the cubic hypersurfaces with
the K3-surfaces is that the moduli of K3-surfaces are much better
understood than those of cubic hypersurfaces them-selves. In
particular, in the case of K3-surfaces one disposes the
surjectivity of the period map. One can expect that a similar
phenomenon holds for the cubic hypersurfaces of dimension $4$.
Would such a result be available, the proofs of the classification
developed in this paper could mimic the proofs of the
corresponding results for K3-surfaces, which would be more simple
and more natural. (The injectivity of the period map and the
Torelli theorem established by C.~Voisin \cite{V} are not
sufficient for our purpose.)
\endrk

\heading Appendix. The list of the eigenlattices of the real
$K3$-involutions.
\endheading

Our list of the eigenlattices  $L_\pm(c)$ of the real
K3-involutions $c\:L\to L$ splits into $3$ tables. In the first
column of each table, K3-vertices $[c]$ are characterized by the
real locus $S_p\+qS$, or $2S_1$, or $\oo$, of the corresponding
K3-surfaces, as explained in 3.4 (cf. also the caption for Figure
1).
 The first two tables describe the so called {\it
 principal series} of lattices $L_\pm$, characterizing the K3-vertices $[S_p\+qS]$
(on Figure 1, these vertices have integer coordinates). The third
table covers the remaining cases and describe $L_\pm$ for
K3-vertices of type $[S_p\+qS]_I$, and in addition $[\oo]$ and
$[2S_1]$ (on Figure 1, these vertices are presented by black
squares which are slightly shifted above the points matching their
values of $(d,r)$). In the first table, we list the lattices $L_+$
fixing in each row the value of $q$ in $[S_p\+qS]$, while in the
second table we list the lattices $L_-$ and fix $p$.

By $\la\pm2\ra$ we mean the lattice of rank $1$ generated by an
element with square $\pm2$, and by $U$ the even unimodular lattice
of rank $2$. Standard notation $E_7$, $E_8$, $D_4$ is used for the
corresponding negative definite root lattices. By $L(2)$ (in our
tables $L=U$ or $E_8$) we denote the lattice whose matrix is the
matrix of $L$ multiplied by $2$. A coefficient, $k$, before a
lattice stands for the direct sum of $k$ copies, for instance,
$2U=U\oplus U$. The direct sum $s\la2\ra\oplus t\la-2\ra$ in the
first two tables is called {\it the diagonal component} of the
corresponding eigenlattice $L_\pm$.

\head Table 1. Principal series of $L_+$ \endhead
 $$\boxed {\matrix
\format\c&\quad\l&\l\\
&\text{diagonal}&\phantom{aa}\text{non-diagonal}\\
K3\text{-vertex}&\text{component}&\phantom{aa}\text{component}\\
\text{---------------}&\text{---------------}&\text{---------------------}\\
[S_{10-t}] &\la2\ra+t\la-2\ra&\\
[S_{10-t}\dsum S] &\phantom{\la2\ra+\,\,}t\la-2\ra&\quad+U\\
[S_{7-t}\dsum 2S] &\phantom{\la2\ra+\,\,}t\la-2\ra&\quad+U+D_4\\
[S_{6-t}\dsum 3S] &\la2\ra+t\la-2\ra&\quad+E_7\\
[S_{6-t}\dsum 4S] &\la2\ra+t\la-2\ra&\quad+E_8\\
[S_{6-t}\dsum 5S] &\phantom{\la2\ra+\,\,}t\la-2\ra&\quad+U+E_8\\
[S_{3-t}\dsum 6S] &\phantom{\la2\ra+\,\,}t\la-2\ra&\quad+U+D_4+E_8\\
[S_{2-t}\dsum 7S] &\la2\ra+t\la-2\ra&\quad+E_7+E_8\\
[S_{2-t}\dsum 8S] &\la2\ra+t\la-2\ra&\quad+2E_8\\
[S_{2-t}\dsum 9S] &\phantom{\la2\ra+\,\,}t\la-2\ra&\quad+U+2E_8\\
\endmatrix}$$

\head Table 2. Principal series of $L_-$ \endhead
$$
\boxed{ \matrix
\format\c&\quad\l&\l\\
&\text{diagonal}&\phantom{aa}\text{non-diagonal}\\
\text{K3-vertex}&\text{component}&\phantom{aa}\text{component}\\
\text{-----------------}&\text{-----------------}&\text{---------------------}\\
[(10-t)S] & 2\la2\ra+t\la-2\ra&\\
[S_1\dsum(9-t)S] &\phantom{2}\la2\ra+t\la-2\ra&\quad+U\\
[S_2\dsum(9-t)S] &\phantom{2\la2\ra+\,\,}t\la-2\ra&\quad  +2U\\
[S_3\dsum(6-t)S]  &\phantom{2\la2\ra+\,\,}t\la-2\ra&\quad  +2U+D_4\\
[S_4\dsum(5-t)S] & 2\la2\ra+t\la-2\ra&\quad  +E_8\\
[S_5\dsum(5-t)S] &\phantom{2}\la2\ra+t\la-2\ra&\quad  +U+E_8\\
[S_6\dsum(5-t)S]  &\phantom{2\la2\ra+\,\,}t\la-2\ra&\quad  +2U+E_8\\
[S_7\dsum(2-t)S]  &\phantom{2\la2\ra+\,\,}t\la-2\ra&\quad  +2U+D_4+E_8\\
[S_8\dsum(1-t)S] & 2\la2\ra+t\la-2\ra&\quad  +2E_8\\
[S_9\dsum(1-t)S] &\phantom{2}\la2\ra+t\la-2\ra&\quad  +U+2E_8\\
[S_{10}\dsum(1-t)S]  &\phantom{2\la2\ra+\,\,}t\la-2\ra&\quad  +2U+2E_8\\
\endmatrix
}$$

\head Table 3. Other types of $L_\pm $ \endhead
$$\boxed{\matrix
\text{K3-vertex}& L_+ &L_-\\
\text{------------------}&\text{---------------------------}&\text{---------------------------}\\
 [8S]_I &U+2D_4+E_8&2U(2)\\
 [S_1\dsum8S]_I &U(2)+2E_8&U+U(2)\\
 [S_1\dsum 4S]_I &U+3D_4&2U(2)+D_4\\
 [S_2\dsum 5S]_I &U(2)+D_4+E_8&U+U(2)+D_4\\
 [S_3\dsum 2S]_I &U(2)+2D_4&U+U(2)+2D_4\\
 [S_4\dsum 3S]_I &U+2D_4&2U+2D_4\\
 [S_5\dsum 4S]_I &U(2)+E_8&U+U(2)+E_8\\
 [S_6\dsum S]_I &U(2)+D_4&U+U(2)+D_4+E_8\\
 [S_9]_I &2U&U+U(2)+2E_8\\
 [2S_1] &U+E_8(2)&2U+E_8(2)\\
 [\oo] &U(2)+E_8(2)&U+U(2)+E_8(2)
\endmatrix}$$

\enddocument